\input amstex
\documentstyle{amsppt}
\magnification=1200

\voffset=-1.5truecm

\TagsOnRight
\NoRunningHeads
\NoBlackBoxes
\loadbold

\define\Ga{\Gamma}
\define\de{\delta}
\define\De{\Delta}
\define\la{\lambda}
\define\La{\Lambda}
\define\om{\omega}
\define\Om{\Omega}
\define\si{\sigma}
\define\ep{\varepsilon}

\define\Alg{\Bbb A}
\define\Ae{\Alg^{\operatorname{ext}}}
\define\C{\Bbb C}
\define\N{\Bbb N}
\define\R{\Bbb R}
\define\Y{\Bbb Y}
\define\Z{\Bbb Z}

\define\ff{\boldkey f}
\define\bh{\boldkey h}
\define\bp{\boldkey p}
\define\tf{\widetilde f}
\define\th{\widetilde h}
\define\tp{\widetilde p}
\define\pp{p^{\#}}
\define\fd{\downarrow}
\define\tA{\widetilde A}
\define\ta{\widetilde a}
\define\hla{\widehat\lambda}

\define\const{\operatorname{const}}
\define\sc{\mu_{\operatorname{s-c}}}
\define\wt{\operatorname{wt}}
\define\Res{\operatorname{Res}}
\define\Span{\operatorname{span}}
\define\sgn{\operatorname{sgn}}
\define\Gauss{\operatorname{Gauss}}
\define\Leb{\operatorname{Leb}}
\define\inv{\operatorname{inv}}

\define\Dz{\Cal D^0}
\define\Mz{\Cal M^0}
\define\m#1{\langle#1\rangle_n}
\define\dc{\overset d\to\longrightarrow}
\define\tht{\thetag}

\topmatter

\title Kerov's central limit theorem for the Plancherel measure on
Young diagrams \endtitle

\author Vladimir Ivanov and Grigori Olshanski \endauthor

\dedicatory In memory of Sergei Kerov (1946--2000) \enddedicatory

\thanks
In: S.~Fomin, editor. Symmetric Functions 2001: Surveys of Developments and
Perspectives (NATO Science Series II. Mathematics, Physics and Chemistry. Vol.
74), Kluwer, 2002, pp. 93--151.
\endgraf
Vladimir Ivanov: Chair of Higher Algebra,
Department of Mathematics and Mechanics, Moscow State
University, Vorob'evy Gory, GZ, Moscow 119992, GSP-2, Russia. E-mail: {\tt
vivanov\@mccme.ru} \endgraf Grigori Olshanski: Dobrushin Mathematics
Laboratory, Institute for Information Transmission Problems, Bolshoy Karetny
19, Moscow 101447, GSP-4, Russia. E-mail: {\tt olsh\@online.ru}
\endthanks

\abstract Consider random Young diagrams with fixed number $n$ of
boxes, distributed according to the Plancherel measure $M_n$. That
is, the weight $M_n(\la)$ of a diagram $\la$ equals $\dim^2\la/n!$,
where $\dim\la$ denotes the dimension of the irreducible representation
of the symmetric group $\frak S_n$ indexed by $\la$. As $n\to\infty$,
the boundary of the (appropriately rescaled) random shape $\la$
concentrates near a curve $\Om$ (Logan--Shepp 1977, Vershik--Kerov
1977). In 1993, Kerov announced a remarkable theorem describing
Gaussian fluctuations around the limit shape $\Om$. Here we propose a
reconstruction of his proof. It is largely based on Kerov's
unpublished work notes, 1999.
\endabstract

\toc
\widestnumber\head{\S10.}
\head \S0. Introduction \endhead
\head \S1. The algebra of polynomial functions on the set of Young
diagrams  \endhead
\head \S2. Continual diagrams and their moments \endhead
\head \S3. The elements $\pp_k$ \endhead
\head \S4. The basis $\{\pp_\rho\}$ and filtrations in $\Alg$ \endhead
\head \S5. The Plancherel measure and the law of large numbers \endhead
\head \S6. The central limit theorem for characters \endhead
\head \S7. The central limit theorem for Young diagrams \endhead
\head \S8. The central limit theorem for transition measures of Young
diagrams \endhead
\head \S9. Discussion \endhead
\head \S10. Free cumulants and Biane's theorem \endhead
\head{} References \endhead
\endtoc

\endtopmatter

\document

\head \S0. Introduction \endhead
\subhead Main result \endsubhead
Let $\Y_n$ denote the set of partitions of $n$ ($n=1,2,\dots$). We
identify partitions and Young diagrams, so that elements of $\Y_n$
become Young diagrams with $n$ boxes. We view each $\la\in\Y_n$ as a
plane shape, of area $n$, inside the first quadrant $\R_+^2$, with
coordinates $r,s$ (the row and column coordinates). In new
coordinates $x=s-r$, $y=r+s$, the boundary $\partial\la$ of the shape
$\la\subset\R_+^2$ may be viewed as the graph of a continuous
piece--wise linear function, which we denote as $y=\la(x)$. Note that
$\la'(x)=\pm1$, and $\la(x)$ coincides with $|x|$ for sufficiently
large values of $|x|$. The area of the shape $|x|\le y\le\la(x)$
equals $2n$.

Further, we equip the finite set $\Y_n$ with a probability measure $M_n$
called the {\it Plancherel measure.\/} The measure $M_n$ has
important representation theoretic and combinatorial interpretations.
By definition, the weight $M_n(\la)$ assigned to a diagram
$\la\in\Y_n$ equals $\dim^2\la/n!$, where $\dim\la$ is the dimension
of the irreducible representation (of the symmetric group $\frak
S_n$) indexed by $\la$. Equivalently, $\dim\la$ is the number of
standard tableaux of shape $\la$.

Viewing $\la$'s as points of the probability space $(\Y_n\,,M_n)$, we
view $\la(\,\cdot\,)$'s as {\it random\/} functions, and we
aim to describe their asymptotics as $n\to\infty$. Informally, the
main result can be stated as follows:
$$
\frac1{\sqrt n}\,\la(\sqrt n\,x)\quad\sim\quad
\Om(x)+\frac2{\sqrt n}\,\De(x), \qquad n\to\infty, \tag0.1
$$
where $\Om(x)$ is a certain fixed curve and $\De(x)$ is a generalized
Gaussian process on the interval $[-2,2]$ ($\Om$ and $\De$ are
specified below).

The left--hand side (denoted as
$\bar\la(x)$ in the sequel) is a {\it rescaled version\/} of the
function $y=\la(x)$. The graph of $\bar\la(\,\cdot\,)$ is obtained
from that of $\la(\,\cdot\,)$ by shrinking both the $x$--axis and the
$y$--axis in $\sqrt n$ times. The purpose of this procedure is to put
the random ensembles with different $n$'s on the same scale (note that the
area of the shape $|x|\le y\le\bar\la(x)$ equals 2 for any $n$).

The first term in the right--hand side of \tht{0.1}
corresponds to the {\it law of large numbers.\/} It follows from
\tht{0.1} that in the large $n$ limit, the random scaled polygonal
lines $y=\bar\la(x)$ concentrate near the fixed curve $y=\Om(x)$.
In the initial scale, this means that, for large $n$, the ``typical''
functions $y=\la(x)$ look like the function
$y=\sqrt n\,\Om(\frac1{\sqrt n} x)$.

The second term in the right--hand side of \tht{0.1} governs the {\it
fluctuations\/} of the random functions $\bar\la$ around the curve
$\Om$, which corresponds to the {\it central limit theorem.} We see that the
fluctuations of the scaled functions are of order $\frac1{\sqrt n}$.
That is, in the initial picture for the shape
$\la\subset\R_+^2$, the random fluctuations of the boundary line
$\partial\la$ need no scaling {\it along the main diagonal\/},
as $n\to\infty$.

\subhead Description of $\Om$ and $\De$ \endsubhead
The function $\Om(x)$ is given by two different expressions depending on
whether $x$ is in the interval $[-2,2]$:
$$
\Om(x)=\cases \frac2\pi(x\arcsin\tfrac x2+\sqrt{4-x^2}), & |x|\le2,\\
|x|, & |x|\ge2.
\endcases \tag0.2
$$
Note that $\Om'(x)=\frac2\pi\arcsin\frac x2$ inside $[-2,2]$.
The critical points $\pm2$ have an important meaning: for
``typical'' (with respect to the Plancherel measure) diagrams
$\la\in\Y_n$, the length of the first row and of the first column is
approximately $2\sqrt n$. This claim, which seems plausible from
\tht{0.1} and \tht{0.2} can be substantially refined, see \cite{BDJ},
\cite{AD}.

The Gaussian process $\De(x)$ can be defined by a random trigonometric
series. Let $\xi_2,\xi_3,\dots$ be independent standard real Gaussian
random variables (each $\xi_k$ has mean 0 and variance 1), and set
$x=2\cos\theta$, where $0\le\theta\le\pi$. Then
$$
\De(x)=\De(2\cos\theta)
=\frac1\pi\,\sum_{k=2}^\infty\frac{\xi_k}{\sqrt k}\,\sin(k\theta).
\tag0.3
$$

For any smooth test function $\varphi$ on $\R$, the smoothed series
$$
\frac1{\pi}\,\sum_{k=2}^\infty \frac{\xi_k}{\sqrt k}\,
\int_{-2}^2 \sin(k\theta)\varphi(x)dx, \qquad
\theta=\arccos(x/2) \tag0.4
$$
converges and is a Gaussian random variable. In this way we get a
Gaussian measure on the space of distributions with support on
$[-2,2]$, or a generalized Gaussian process. Its trajectories are not
ordinary functions but generalized functions.

\subhead History of the result \endsubhead
The law of large numbers (the concentration near the curve $\Om$) was
independently obtained by Logan and Shepp \cite{LoS} and by Vershik and
Kerov \cite{VeK1}. Their papers appeared in 1977. Later, in 1985,
Vershik and Kerov published a detailed version of their work,
\cite{VeK3}, containing stronger results. In \cite{LoS} and
\cite{VeK3}, the question about the second term of the asymptotics,
corresponding to the central limit theorem, was posed.

Such a theorem was obtained by Kerov and announced in his short
note \cite{Ke1}, 1993. There Kerov also outlined the scheme of the proof.
The note \cite{Ke1} contained a number of fruitful ideas, one of which
(introduction of ``good'' coordinates in the set of Young diagrams)
was largely developed in the joint note by Kerov and Olshanski
\cite{KO}, 1994.

For an intermediate result of \cite{Ke1}, which is of independent
interest, an elegant proof was suggested by Hora \cite{Ho}, 1998.
Note that Hora's approach differs from that of Kerov.

A few years ago we started to persuade
Kerov to write a detailed exposition of his central limit theorem.
Our discussions resulted first in the joint paper by Ivanov and
Kerov \cite{IK}, 1999, which clarified and developed one of the steps
of Kerov's proof. \footnote{It concerns the stable structure constants
for convolution of conjugacy classes in symmetric groups. This topic
was also discussed in \cite{KO}.}

Then Kerov found a simpler derivation of the theorem,
which also made apparent that the subject is connected with the
concept of free cumulants and a theorem due to Biane \cite{Bi1}. In
the end of 1999 Kerov sent us two short work notes with a description
of the new approach. About the same time he gave a talk on this
subject at Vershik's seminar in St.~Petersburg. He also started
writing a detailed paper on this subject but had time only to finish the
preliminary section.

In the present paper we give a detailed exposition of Kerov's central
limit theorem. Our aim was to reconstruct from his notes the ``new
approach'' of 1999. \footnote{It should be pointed out that the
``old approach'' of 1993 is also correct: we were able to
directly check all the claims of \cite{Ke1}.} This was not easy: for a
long time we could not understand the meaning of some claims stated
too briefly, but finally the picture became clear.
However, we cannot be sure that we succeeded to completely fathom
Kerov's intention, and there is no doubt that his own exposition
would be quite different.

\subhead Links with random matrices \endsubhead
Recently it was discovered that the limit distribution of
a finitely many (properly scaled) first rows of the random Plancherel
Young diagram $\la\in\Y_n$, as $n\to\infty$, coincides with the limit
distribution of the same number of (properly scaled) largest
eigenvalues of the random Hermitian $N\times N$ matrix taken from the
Gaussian Unitary Ensemble, as $N\to\infty$. See \cite{BDJ} and
subsequent papers \cite{Ok}, \cite{BOO}, \cite{Jo3}, \cite{Jo4},
\cite{BDR}. It turns out that the striking similarity between these
two random ensembles holds not only ``at the edge'' (as is shown in
these works) but also on the level of global fluctuations, which is
the subject of the present paper. For spectra of random matrices, the
limit behavior of global fluctuations was first studied in \cite{DS},
where a central limit theorem was obtained. Further results in this
direction were obtained in \cite{Jo1}, \cite{Jo2}, \cite{DE}. The
generalized Gaussian processes that emerge in these works are very
close to our process, we discuss this topic in \S9.

\subhead Techniques \endsubhead
Although the main result is stated in probabilistic terms, the
techniques of the paper are essentially algebraic and combinatorial,
the probabilistic part being reduced to a few elementary facts. The
work is based on the choice of convenient ``coordinate systems'' for
Young diagrams (there are several ones) and on the choice of an
appropriate algebra $\Alg$ of ``observables''. Elements of $\Alg$ are
functions on the set $\Y$ of all Young diagrams. They are given by
polynomial expressions in each of the ``coordinate systems''. For this
reason we call $\Alg$ the algebra of polynomial functions on $\Y$. We
examine several different bases in $\Alg$. One of them (denoted as
$\{\pp_\rho\}$) is related to the character table of the symmetric
groups; this basis is well adapted to evaluating expectations with
respect to the Plancherel measures $M_n$. Another basis has
geometric significance; this basis is formed by monomials in
$\tp_2,\tp_3,\dots$, a system of generators of $\Alg$, which are
essentially the moments of $\la(x)$. One more basis (formed by
monomials in generators $p_1,p_2,\ldots\in\Alg$) plays an intermediate
role. A major part of our work consists in studying the transitions
between various bases. This finally makes it possible to isolate a
good system of generators in $\Alg$ that directly describe the Gaussian
fluctuations.

\subhead Organization of the paper \endsubhead
In \S1, we introduce the algebra $\Alg$ and a system $p_1,p_2,\dots$
of its generators. We show that elements of $\Alg$ are both shifted
symmetric functions in the row coordinates $\la_1,\la_2,\dots$ of a
Young diagram $\la\in\Y$, and supersymmetric functions in the
(modified) Frobenius coordinates of $\la$. This fact was first
pointed out in \cite{KO}.

In \S2, we introduce the necessary geometric setting for visualizing
fluctuations of Young diagrams. We embed $\Y_n$ into the larger set $\Dz$ of
``continual diagrams''. We introduce the generators $\tp_2,\tp_3,\ldots\in\Alg$
and the ``weight grading'' of the algebra $\Alg$, which is well adapted to the
operation of rescaling diagrams.

In \S3, we examine one more system of generators in $\Alg$, denoted
as $\pp_1,\pp_2,\dots$\,. These are character values on cycles in
symmetric groups. We study the transitions between all three
systems of generators. Here our tools are a suitably elaborated
classical formula (due to Frobenius) for the value of a symmetric
group character on the $k$-cycle, and Lagrange's inversion formula.

In \S4, we introduce the basis $\{\pp_\rho\}$ in $\Alg$ and study a
family of filtrations in $\Alg$, which are defined in terms of
this basis. Here we follow the paper \cite{IK}. We essentially need two
different filtrations. Their purpose is to single out main terms of
asymptotics in different regimes. One filtration is responsible for
the ``law of large numbers'' while another serves the
``central limit theorem''.

In \S5, we start the study of the Plancherel measures $M_n$. We
introduce the sequence of expectation functionals $\m{\,\cdot\,}$ on
$\Alg$ that corresponds to the sequence $\{M_n\}$, and we remark that
$\m{\,\cdot\,}$ becomes very simple in the basis $\{\pp_\rho\}$. Then
we prove the main result of the section --- the law of large
numbers, or convergence to the curve $\Om$. Although the
central limit theorem, established in \S7, contains the law of large
numbers, we prefer to prove it independently, because this can be
done in a rather simple way. It is interesting to compare our
simple algebraic argument with the analytic approach of the pioneer
works \cite{LoS}, \cite{VeK1}, \cite{VeK3}.

In \S6, we examine the random variables
${\pp_2}^{(n)},{\pp_3}^{(n)},\dots$, where ${\pp_k}^{(n)}$ stands for the
restriction of the function $\pp_k\in\Alg$ to the finite probability
space $(\Y_n,M_n)$. We show that, as $n\to\infty$, the variables
${\pp_k}^{(n)}$, suitably scaled, are asymptotically independent
Gaussians. This result is the first version of the central limit
theorem. Its proof relies on the method of \cite{IK}. A different
proof has been given by Hora \cite{Ho}.

In \S7, we obtain our main result: a description of the Gaussian
fluctuations around the limit curve $\Om$. It is derived from the
central limit theorem for the generators $\pp_k$ mentioned above. The
proof is based on a formula that gives the highest term of the
polynomial expressing $\pp_k$ through the (centered and scaled
versions of) the generators $\tp_j$. Here ``highest term'' refers to
an appropriate filtration of the algebra $\Alg$, which we call
Kerov's filtration.

In \S8, we get one more version of the central limit theorem.
According to Vershik--Kerov's theory, to any Young diagram $\la$
we attach a probability measure on $\R$ (say, $\mu_\la$), supported
by a finite set. Viewing
$\la\in\Y_n$ as the random element of the probability space
$(\Y_n,M_n)$, we turn $\mu_\la$ into a random measure. For these random
measures we prove an asymptotic formula similar to \tht{0.1}, where, instead
of the limit curve $\Om$, we have the semi--circle distribution, and
$\De(x)$ is replaced by another generalized Gaussian process.
We do not know if Kerov was aware of this result. However, it
perfectly fits in the philosophy of his works.

In \S9, we give comments to the results of \S\S7--8 and compare them with
the central limit theorem for random matrices.

In \S10, we show that the highest terms of the elements $\pp_\rho$ in
the ``weight grading'' are closely related to
the free cumulants. As an application, we get a simple proof of Biane's
asymptotic formula for character values of large symmetric groups,
\cite{Bi1}.

\subhead Acknowledgment \endsubhead
One of the authors (G.~O.) is deeply grateful to Persi Diaconis for
discussions and an important critical remark, which was taken into
account in the final version of the paper.

\head \S1. The algebra of polynomial functions on the set of Young
diagrams \endhead

Recall first the basic definitions and notation related to partitions
and Young diagrams, see \cite{Ma}.

A {\it partition \/} is an infinite sequence
$\la=(\la_1,\la_2,\dots)$ of nonnegative integers such that
$\la_1\ge\la_2\ge\dots$ and the number of nonzero $\la_i$'s is
finite. The sum $\la_1+\la_2+\dots$ is denoted by $|\la|$, and
usually we set $|\la|=n$.

As in \cite{Ma}, we assign to a partition a Young diagram, which is
denoted by the same symbol. We identify partitions and Young
diagrams, and we denote by $\Y$ the set of all Young diagrams. The
conjugation involution of $\Y$ (transposition of rows and columns of
a diagram) is denoted as $\la\mapsto\la'$.

There is another presentation of Young diagrams, the {\it Frobenius
notation.\/} We shall use its modification due to Vershik and Kerov
\cite{VeK2}:
$$
\la=(a_1,\dots,a_d\mid b_1,\dots,b_d).
$$
Here $d=d(\la)$ is the length of the main diagonal in $\la$,
$$
d(\la)=\{i\mid \la_i\ge i\}=\{j\mid \la'_j\ge j\},
$$
and
$$
a_i=\la_i-i+\tfrac12, \quad b_i=\la'_i-i+\tfrac12,
\qquad i=1,\dots, d(\la). \tag1.1
$$

The numbers $a_i, b_i$ are called the {\it modified Frobenius
coordinates\/} of $\la$. Both $a_1,\dots,a_d$ and $b_1,\dots,b_d$ are
strictly decreasing positive proper half--integers, i.e.,
numbers from $\{\tfrac12,\tfrac32,\tfrac52,\dots\}$, such that $\sum
(a_i+b_i)=|\la|$. Note that the Vershik--Kerov definition \tht{1.1}
differs from the classical definition of the Frobenius
coordinates, which does not involve one--halves (see \cite{Ma, p. 3}).
However, these one--halves play an important role in what follows.

Note that the conjugation involution $\la\mapsto\la'$ has a very
simple description in terms of the Frobenius coordinates:
$$
(a_1,\dots,a_d\mid b_1,\dots,b_d)'=
(b_1,\dots,b_d\mid a_1,\dots,a_d). \tag1.2
$$

One more useful presentation of Young diagrams, due to Kerov,
will be given in \S2.

Set
$$
\Z'=\Z+\tfrac12
=\{\dots,-\tfrac32,-\tfrac12,\tfrac12,\tfrac32,\dots\},\quad
\Z'_+=\{\tfrac12,\tfrac32,\dots\},\quad
\Z'_-=\{\ldots,-\tfrac32,-\tfrac12\}.
$$

Given $\la\in\Y$, set
$$
l_i=\la_i-i+\tfrac12\in\Z', \quad i=1,2,\dots,
$$
and note that $l_1>l_2>\dots$. We assign to $\la$ the infinite subset
$\Cal L(\la)=\{l_1,l_2,\dots\}\subset\Z'$.

The following claim is a version of the classical Frobenius lemma,
see \cite{Ma, ch. I, (1.7) and Example 1.15 (a)}.

\proclaim{Proposition 1.1} Let $\la\in\Y$ be arbitrary.

{\rm(i)} We have $\Z'=\Cal L(\la)\sqcup(-\Cal L(\la'))$. I.e., $\Cal
L(\la)\cap(-\Cal L(\la'))=\varnothing$ and
$\Cal L(\la)\cup(-\Cal L(\la'))=\Z'$.

{\rm(ii)} In the notation \tht{1.1},
$$
\Cal L(\la)\cap\Z'_+=\{a_1,\dots,a_d\},\quad
\Z'_-\setminus\Cal L(\la)=-(\Cal L(\la')\cap\Z'_+) =\{-b_1,\dots,-b_d\}).
$$
\endproclaim

\demo{Proof} (i) We represent $\la$ as a plane shape in the quarter
plane $\R_+^2$. Let $(r,s)$ be the coordinates in $\R_+^2$. Here the
rows of $\la$ are counted along the first coordinate $r$, directed
downwards, while the columns are counted along the second coordinate
$s$, directed to the right. Denote by $\partial\la$ the doubly
infinite polygonal line which first goes upwards along the $r$-axis, next
goes along the boundary line separating $\la$ from
its complement in $\R_+^2$, and then goes to the right along the
$s$--axis. For any $a\in\Z'$, the diagonal line $s-r=a$ intersects
$\partial\la$ at the midpoint of a certain segment, which is either vertical
or horizontal. According to these two possibilities $a$ is either in
$\Cal L(\la)$ or in $-\Cal L(\la')$. This proves (i).

(ii) By the very definition \tht{1.1}, the numbers $a\in\Z'_+$ such that
the diagonal $s-r=a$ meets a vertical boundary segment are exactly
the numbers $a_1,\dots,a_d$. Likewise, the numbers $-b\in\Z'_-$ such
that the diagonal $s-r=-b$ meets a horizontal boundary segment are
exactly the numbers $-b_1,\dots,-b_d$. This proves (ii). \qed
\enddemo

For any $\la\in\Y$ we set
$$
\Phi(z;\la)=\prod_{i=1}^\infty\frac{z+i-\tfrac12}{z-\la_i+i-\tfrac12}\,,
\qquad z\in\C.
$$
The product is actually finite, because $\la_i=0$ when $i$ is large
enough. Therefore, $\Phi(z;\la)$ is a rational function in $z$. We
view it as a generating function of $\la$.

\proclaim{Proposition 1.2} In the notation \tht{1.1}, we have
$$
\Phi(z;\la)=\prod_{i=1}^d\frac{z+b_i}{z-a_i}\,, \tag1.3
$$
which is the presentation of $\Phi(z;\la)$ as an incontractible
fraction.
\endproclaim

\demo{Proof} The equality \tht{1.3} follows from Proposition 1.1. This
is an incontractible fraction, because the numbers
$a_1,\dots,a_d,-b_1,\dots,-b_d$ are pairwise distinct. \qed
\enddemo

Another proof of the equality \tht{1.3} is given in \cite{ORV}, it
follows an idea from \cite{KO}.

As the first corollary of \tht{1.3} note the relation
$$
\Phi(z;\la')=1/\Phi(-z;\la).
$$

Remark that $\Phi(z;\la)=1+O(\tfrac 1z)$ near $z=\infty$, hence both
$\Phi(z;\la)$ and $\ln\Phi(z;\la)$ can be expanded in a power series
in $z^{-1},z^{-2},\dots$ about $z=\infty$.

\example{Definition 1.3} The {\it algebra of polynomial functions on
the set $\Y$,\/} denoted as $\Alg$, is generated over $\R$ by the coefficients
of the above expansion of $\Phi(z;\la)$ or, equivalently, of
$\ln\Phi(z;\la)$. We also assume that $\Alg$ contains 1.
\endexample

\proclaim{Proposition 1.4} We have
$$
\ln\Phi(z;\la)=\sum_{k=1}^\infty\frac{p_k(\la)}k\,z^{-k}\,,
$$
where
$$
\align
p_k(\la)&=\sum_{i=1}^\infty [(\la_i-i+\tfrac12)^k-(-i+\tfrac12)^k]\tag1.4\\
&=\sum_{i=1}^{d(\la)}[a_i^k-(-b_i)^k]. \tag1.5
\endalign
$$
\endproclaim

\demo{Proof} Immediate from Proposition 1.2.\qed
\enddemo

Thus, $\Alg$ is generated by the functions $p_k(\la)$, $k=1,2,\dots$.

Recall \cite{Ma} that the {\it algebra of symmetric functions,\/}
denoted as $\La$, is the graded algebra defined as the projective limit
(in the category of graded algebras) $\varprojlim\La(n)$, where
$\La(n)$ denotes the algebra of symmetric polynomials in $n$
variables. As the base field we take $\R$. The morphism
$\La(n)\to\La(n-1)$, which is employed in the projective limit
transition, is defined as specializing the $n$th variable to 0. Let
$\{\bh_k\}_{k=1,2,\dots}$ and $\{\bp_k\}_{k=1,2,\dots}$ denote the
complete homogeneous symmetric functions and the Newton power sums,
respectively. Each of these two families is a system of homogeneous,
algebraically independent generators of $\La$, $\deg\bh_k=\deg\bp_k=k$.
Recall the basic relation:
$$
1+\sum_{k=1}^\infty \bh_kt^k
=\exp\sum_{k=1}^\infty \frac{\bp_k}k\,t^k\,.
$$

\proclaim{Proposition 1.5} The generators $p_k\in\Alg$ are algebraically
independent, so that $\Alg$ is isomorphic to $\R[p_1,p_2,\dots]$.
\endproclaim

\demo{Proof} Fix an arbitrary $N=1,2,\dots$. Assume that $f$ is a
polynomial in $N$ variables such that
$f(p_1,\dots,p_N)=0$, and show that $f=0$. Let $\bar f$ denote the top
homogeneous component of $f$ counted with the understanding that the
degree of the $i$th variable equals $i$; it suffices to show that
$\bar f=0$.

Let $\la$ range over the set of partitions of length $\le N$. Fix an
arbitrary vector $x\in\R^N$ with nonnegative weakly decreasing
coordinates and set $\la=\la(A)=([Ax]_i)_{i=1,\dots,N}$\,, where $A$ is
a large integer. Letting $A\to\infty$ in the equality
\linebreak
$f(p_1(\la(A)),\dots,p_N(\la(A)))=0$ we get $\bar
f(\bp_1(x),\dots,\bp_N(x))=0$. Since the first $N$ Newton power sums
specialized in $N$ variables are algebraically independent, we
conclude that $\bar f=0$. \qed
\enddemo

\example{Definition 1.6} Setting $\La\ni\bp_k\mapsto p_k\in\Alg$ and
taking into account Proposition 1.5 we get an algebra isomorphism
$\La\to\Alg$. We call it the {\it canonical isomorphism.\/} We call the
grading in $\Alg$, inherited from that of $\La$, the {\it canonical
grading\/} of $\Alg$.
\endexample

Later on, in Definition 2.9, we shall define quite a different
grading in $\Alg$.

In terms of generating series, the canonical isomorphism $\La\to\Alg$
takes the form
$$
H(t):=1+\sum_{k=1}^\infty\bh_kt^k\, \mapsto \,
\Phi(t^{-1};\,\cdot\,).
$$

Formula \tht{1.5} means that the functions $p_k(\la)$ are {\it super\/}
power sums in $a_i$'s and $b_i$'s, see \cite{Ma, Example I.3.23},
\cite{VeK2}, \cite{KO}, \cite{ORV}. Thus, one
can say that under the canonical isomorphism of Definition 1.6, the
algebra $\Alg$ is identified with the {\it algebra of supersymmetric
functions in the modified Frobenius coordinates\/} of a Young diagram.

Next, we shall give a similar interpretation of formula \tht{1.4}.
Recall \cite{OO} that the {\it algebra of shifted symmetric
functions,\/} denoted as $\La^*$, is the filtered algebra defined as
the projective limit (in the category of filtered algebras)
$\varprojlim\La^*(n)$, where $\La(n)^*$ consists of
those polynomials in $n$ variables $x_1,\dots,x_n$, which become
symmetric in new variables $y_i=x_i-i+\const$ (the choice of the
constant here is irrelevant). The base field is again $\R$, the
filtration is taken with respect to the total degree of a polynomial,
and the morphism $\La^*(n)\to\La^*(n-1)$ is defined as above, i.e.,
as specializing $x_n=0$. The graded algebra associated to the
filtered algebra $\La^*$ is canonically isomorphic to $\La$.
The algebra $\La^*$ is generated by the algebraically independent
system $\{\bp^*_k\}_{k=1,2,\dots}$, where
$$
\bp^*_k(x_1,x_2,\dots)
=\sum_{k=1}^\infty[(x_i-i+\tfrac12)^k-(-i+\tfrac12)^k],
\quad k=1,2,\dots,
$$
are certain shifted analogs of the Newton power sums. See \cite{OO}
for more detail (note that the above definition of the elements
$\bp^*_k$ slightly differs from that given in \cite{OO}). See also
\cite{EO}.

By analogy with Definition 1.6, we define an algebra isomorphism
$\La^*\to\Alg$ by setting $\bp^*_k\mapsto p_k$, $k=1,2,\dots$. Note
that it preserves the filtration. Then formula \tht{1.4} makes it
possible to say that the algebra $\Alg$ coincides with the {\it algebra of
shifted symmetric functions in the row coordinates\/}
$\la_1,\la_2,\dots$ of a Young diagram $\la$.

\example{Definition 1.7} Define an involutive algebra automorphism
$\inv:\Alg\to\Alg$ by
$$
(\inv(f))(\la)=f(\la'), \qquad f\in\Alg, \quad\la\in\Y.  \tag1.6
$$
\endexample

By virtue of \tht{1.2} and \tht{1.5},
$$
\inv(p_k)=(-1)^{k-1}p_k\,,\qquad k=1,2,\dots\,. \tag1.7
$$
Hence the involution of $\Alg$ is compatible with the canonical involution
of the algebra $\La$ with respect to the isomorphism $\La\to\Alg$
introduced in Definition 1.6.

\head \S2. Continual diagrams and their moments \endhead

\example{Definition 2.1} A {\it continual diagram\/} is a function
$\om(x)$ on $\R$ such that:

(i) $|\om(x_1)-\om(x_2)|\le |x_1-x_2|$ for any $x_1,x_2\in\R$ (the
Lipschitz condition).

(ii) There exists a point $x_0\in\R$, called the {\it center\/} of
$\om$, such that $\om(x)=|x-x_0|$ when $|x|$ is large enough.

The set of all continual diagrams is denoted by $\Cal D$, and the subset
of diagrams with center 0 is denoted by $\Dz$.
\endexample

This definition is due to Kerov (see his papers \cite{Ke2},
\cite{Ke3}, \cite{Ke4}). We shall mainly deal with the set $\Dz$.

To any $\om\in\Cal D$ we assign a function $\si(x)$:
$$
\si(x)=\tfrac12(\om(x)-|x|). \tag2.1
$$
Since $\si(x)$ satisfies the Lipschitz condition (i), its derivative
$\si'(x)$ exists almost everywhere and satisfies $|\si'(x)|\le1$. By
(ii), the function $\si'(x)$ is compactly supported.

If $\om\in\Dz$ then $\si(x)$ is compactly supported, too. For general
$\om\in\Cal D$, we have
$$
\si(x)\equiv-x_0, \quad x\gg0; \qquad
\si(x)\equiv x_0, \quad x\ll0.
$$
This implies that $\om(x)$ is uniquely determined by $\si'(x)$. Even
more, $\om(x)$ is uniquely determined by the second derivative
$\si''(x)$, which is understood in the sense of distribution theory.

Define the functions $\tp_1,\tp_2,\dots$ on $\Cal D$ by setting
$$
\align
\tp_k[\om]&= -k\int_{-\infty}^\infty x^{k-1}\si'(x)dx  \tag2.2\\
&=\int_{-\infty}^\infty x^k\si''(x)dx, \tag2.3
\endalign
$$
where $\om\in\Cal D$, $k=1,2,\dots$.

\proclaim{Proposition 2.2} If $\om\in\Dz$ then $\tp_1[\om]=0$ and
$$
\tp_k[\om]=k(k-1)\int_{-\infty}^\infty x^{k-2}\si(x)dx, \quad
k=2,3,\dots
$$
\endproclaim

\demo{Proof} Recall that $\si(x)$ is finitely supported when
$\om\in\Dz$. This implies the first claim. Further, integrating
\tht{2.2} by parts gives the second claim. \qed
\enddemo

\example{Definition 2.3} Given $\la\in\Y$, we define a piece--wise
linear function
$\la(\,\cdot\,)$ as follows. Let $(r,s)$ and $\partial\la$ be
as in the proof of Proposition 1.1. Then the
graph $y=\la(x)$ describes $\partial\la$ in the coordinates
$x=s-r$, $y=r+s$. The correspondence $\la\mapsto\la(\,\cdot\,)$ yields an
embedding $\Y\hookrightarrow\Dz$.
\endexample

We have $\la'(x)=\pm1$, except finitely many points, which are
exactly the local extrema of the function $\la(x)$. These local
extrema form two interlacing sequences of points
$$
x_1<y_1<x_2<\dots<x_m<y_m<x_{m+1}\,, \tag2.4
$$
where the $x_i$'s are the local minima and the $y_j$'s are the local
maxima of the function $\la(x)$.

\proclaim{Proposition 2.4} We have
$$
x_i\in\Z, \quad y_j\in\Z, \quad \sum x_i-\sum y_j=0. \tag2.5
$$
Conversely, any couple of interlacing sequences \tht{2.4} satisfying
\tht{2.5} comes from a Young diagram $\la$, which is determined
uniquely.
\endproclaim

\demo{Idea of proof} For any couple of interlacing sequences
\tht{2.4}, there exists a unique polygonal line $\om\in \Cal D$, with
center at $x_0=\sum x_i-\sum y_j$ and such that, for the
corresponding function $\si$,
$$
\si''(x)=\sum_{i=1}^{m+1}\delta(x-x_i)-\sum_{j=1}^m \delta(x-y_j)
-\delta(x-x_0).
$$
The line $\om$ represents a Young diagram if and only if $x_0=0$. \qed
\enddemo

The correspondence $\la\mapsto\{x_i\}\cup\{y_j\}$ provides one more
useful system of parameters for Young diagrams.

\proclaim{Proposition 2.5} Let $\la\in\Y$, let
$\la(\,\cdot\,)\in\Dz$ be the corresponding continual diagram, and
consider the local extrema \tht{2.4}. We have
$$
\tp_k[\la(\,\cdot\,)]=\sum_{i=1}^{m+1}x_i^k-\sum_{j=1}^m y_j^k\,,
\quad k=1,2,\dots.
$$
\endproclaim

\demo{Proof} Let $\si$ be associated with $\om=\la(\,\cdot\,)$, as
defined in \tht{2.1}. Then we get
$$
\si''(x)=\sum_{i=1}^{m+1}\delta(x-x_i)-\sum_{j=1}^m \delta(x-y_j)
-\delta(x).
$$
Note that $\int x^k\delta(x)dx=0$ for any $k=1,2,\dots$ and apply \tht{2.3}.
\qed
\enddemo

\proclaim{Proposition 2.6} Let $\la\in\Y$ and let $\{x_i\}\cup\{y_j\}$ be
the local extrema of $\la(\,\cdot\,)$. The following identity holds
$$
\frac{\Phi(z-\tfrac12;\la)}{\Phi(z+\tfrac12;\la)}
=\frac{z\prod_{j=1}^m(z-y_j)}{\prod_{i=1}^{m+1}(z-x_i)}\,. \tag2.6
$$
\endproclaim

\demo{Proof} We shall prove the identity
$$
\Phi(z-\tfrac12;\la)
=\frac{\prod_{i=1}^{m+1}\Ga(z-x_i)}{\Ga(z)\prod_{j=1}^m\Ga(z-y_j)}\,,
\tag2.7
$$
which implies \tht{2.6}. Using \tht{1.3}, we rewrite \tht{2.7} as
$$
\prod_{i=1}^d\frac{z+b_i}{z-a_i}
=\frac{\prod_{i=1}^{m+1}\Ga(z-x_i+\tfrac12)}
{\Ga(z+\tfrac12)\prod_{j=1}^m\Ga(z-y_j+\tfrac12)}\,. \tag2.8
$$
Let $(r,s)$ be the row and column coordinates in the quarter--plane,
see the proof of Proposition 1.1. Draw the diagonal lines $s-r=x_i$
($1\le i\le m+1$) and $s-r=y_j$ ($1\le j\le m$), which divide the
boundary line $\partial\la$ into interlacing vertical and horizontal
pieces.

Assume first that $x_l<0<y_l$ for a certain $l$. Consider an
arbitrary {\it vertical\/} piece of $\partial\la$ which sits {\it
above\/} the main diagonal $s-r=0$. The ends of such a piece lie on
the lines $s-r=y_k$ and $s-r=x_{k+1}$, where $l\le k\le m$. Inside
this piece, the row Frobenius coordinates increase by one and form
the sequence
$$
y_k+\tfrac12,\quad
y_k+\tfrac32,\quad \dots, \quad
x_{k+1}-\tfrac32,\quad
x_{k+1}-\tfrac12\,,
$$
so that the partial product in $\prod_{i=1}^d 1/(z-a_i)$
corresponding to this sequence equals
$$
\frac{\Ga(z-x_{k+1}+\tfrac12)}{\Ga(z-y_k+\tfrac12)}\,.
$$
It follows that
$$
\prod_{i=1}^d\frac1{z-a_i}=
\prod_{k=l}^m\frac{\Ga(z-x_{k+1}+\tfrac12)}{\Ga(z-y_k+\tfrac12)}
=\frac{\prod_{k=l+1}^{m+1}\Ga(z-x_k+\tfrac12)}
{\prod_{k=l}^m\Ga(z-y_k+\tfrac12)}\,. \tag2.9
$$

Next, consider a {\it horizontal\/} piece {\it below\/} the main
diagonal. Such a piece sits between the diagonal lines $s-r=x_k$ and
$s-r=y_k$, where $1\le k\le l-1$. In this piece, the column Frobenius
coordinates make up the sequence
$$
-(x_k+\tfrac12),\quad
-(x_k+\tfrac32),\quad \dots,\quad
-(y_k-\tfrac32), \quad
-(y_k-\tfrac12),
$$
whose contribution to the product $\prod_{i=1}^d(z+b_i)$ equals
$$
\frac{\Ga(z-x_k+\tfrac12)}{\Ga(z-y_k+\tfrac12)}\,.
$$
Therefore, the contribution of all horizontal pieces below the main
diagonal equals
$$
\prod_{k=1}^{l-1}\frac{\Ga(z-x_k+\tfrac12)}{\Ga(z-y_k+\tfrac12)}\,. \tag2.10
$$

Finally, consider the only piece that {\it intersects\/} the main
diagonal. By our assumption, this piece is {\it horizontal\/}, and it
sits between the lines $s-r=x_l$ and $s-r=y_l$. We have to examine
the row Frobenius coordinates inside it. They make up the sequence
$$
-(x_l+\tfrac12), \quad
-(x_l+\tfrac32), \quad \dots, \quad \tfrac12
$$
(recall that $x_l<0<y_l$). The corresponding contribution equals
$$
\frac{\Ga(z-x_l+\tfrac12)}{\Ga(z+\tfrac12)}\,. \tag2.11
$$

Multiplying up \tht{2.10} and \tht{2.11} we get
$$
\gather
\prod_{i=1}^d(z+b_i)=
\prod_{k=1}^{l-1}\frac{\Ga(z-x_k+\tfrac12)}{\Ga(z-y_k+\tfrac12)}
\,\cdot\,\frac{\Ga(z-x_l+\tfrac12)}{\Ga(z+\tfrac12)}\\
=\frac{\prod_{k=1}^l\Ga(z-x_k+\tfrac12)}
{\Ga(z+\tfrac12)\,\prod_{k=1}^{l-1}\Ga(z-y_k+\tfrac12)}\,. \tag2.12
\endgather
$$
Now \tht{2.8} follows from \tht{2.9} and \tht{2.12}.

We have verified \tht{2.6} under the assumption $x_l<0<y_l$, i.e.,
in the case when the main diagonal $s-r=0$ meets the boundary line
$\partial\la$ at an interior point of a horizontal piece. The same
argument works if the intersection of $s-r=0$ with $\partial\la$ is
inside a vertical piece (i.e., $y_l<0<x_{l+1}$ for a certain $l$) or
if the intersection point coincides with a brake of $\partial\la$
(i.e., some $x_l$ or $y_l$ is 0). \qed
\enddemo

Set
$$
\tp_k(\la)=\tp_k[\la(\,\cdot\,)], \qquad
k=1,2,\dots, \quad \la\in\Y,
$$
where the right--hand side is given by \tht{2.2} or, equivalently, by
\tht{2.3}. Note that $\tp_1(\la)\equiv0$.

\proclaim{Proposition 2.7} The functions $\tp_2(\la), \tp_3(\la),\dots$ belong
to the algebra $\Alg$ and are related to the functions
$p_1(\la),p_2(\la)$ by the relations
$$
\tp_k=\sum_{j=0}^{[\frac{k-1}2]}\binom
k{2j+1}\,2^{-2j}\,p_{k-1-2j}\,, \qquad k=2,3,\dots\,. \tag2.13
$$
\endproclaim

\demo{Proof} Formula \tht{2.6} implies that
$$
\ln\Phi(z-\tfrac12;\la)-\ln\Phi(z+\tfrac12;\la)
=\ln\prod_{j=1}^m\left(1-\frac{y_j}z\right)
-\ln\prod_{i=1}^{m+1}\left(1-\frac{x_i}z\right). \tag2.14
$$
By Proposition 1.4, the left--hand side equals
$$
\gather
\sum_{l=1}^\infty \frac{p_l(\la)}l\,
\left(\frac1{(z-\tfrac12)^l}-\frac1{(z+\tfrac12)^l}\right)\\
=\sum_{l=1}^\infty\frac{p_l(\la)}l\,z^{-l}\,
\left(\left(1-\frac1{2z}\right)^{-l}-\left(1+\frac1{2z}\right)^{-l}\right)\\
=\sum_{l=1}^\infty\sum_{j=0}^\infty\frac{l(l+1)\dots(l+2j)}{(2j+1)!}\,
\frac{p_l(\la)}l\, \frac{z^{-(l+2j+1)}}{2^{2j}}\\
=\sum_{l=1}^\infty\sum_{j=0}^\infty\frac{(l+1)\dots(l+2j+1)}{(2j+1)!}\,
\frac{p_l(\la)}{l+2j+1}\, \frac{z^{-(l+2j+1)}}{2^{2j}}\,.
\endgather
$$
Setting $l+2j+1=k$ we rewrite this as
$$
\sum_{k=2}^\infty \sum_{j=0}^{[\frac{k-1}2]}\binom k{2j+1}\,
\frac{p_{k-2j-1}(\la)}{2^{2j}}\,
\frac{z^{-k}}k\,. \tag2.15
$$
By Proposition 2.5, the right--hand side of \tht{2.14} equals
$$
\sum_{k=1}^\infty\tp_k(\la)\,\frac{z^{-k}}k\,. \tag2.16
$$
Comparing the coefficients of $z^{-k}/k$ in \tht{2.15} and \tht{2.16} we get
$\tp_1(\la)\equiv0$ (which we already know) and then \tht{2.13}. \qed
\enddemo

Note that
$$
\inv(\tp_k)=(-1)^k\tp_k\,, \qquad k=2,3,\dots, \tag2.17
$$
where `$\inv$' is the involution introduced in Definition 1.7. Indeed,
\tht{2.17} easily follows from the definition of $\tp_k$ and the
symmetry property $\la'(x)=\la(-x)$. The fact that in the right--hand
side of \tht{2.13}, the subscript varies with step 2 agrees with the
symmetry properties of $p_k$'s and $\tp_k$'s, see \tht{1.7} and
\tht{2.17}.

\proclaim{Corollary 2.8} For any $k=2,3,\dots$
$$
\frac{\tp_k}k=p_{k-1}\,+\,
\langle\text{\rm a linear combination of $p_{k-2},
\dots,p_1$}\rangle.
$$

Conversely, for any $k=1,2,\dots$
$$
p_k=\frac{\tp_{k+1}}{k+1}\, +\,
\langle\text{\rm a linear combination of
$\tp_k,\dots,\tp_2$}\rangle.
$$
\endproclaim
\qed

By Corollary 2.8, the elements $\tp_2,\tp_3,\dots$ are algebraically
independent generators of the algebra $\Alg$:
$$
\Alg=\R[\tp_2,\tp_3,\dots].
$$

\example{Definition 2.9 (cf. \cite{EO})} The {\it weight grading\/}
of the algebra $\Alg$ is defined by setting
$$
\wt(\tp_k)=k, \qquad k=2,3,\dots
$$
Equivalently, the weight grading is the image of the standard grading
of $\La$ under the algebra morphism
$$
\gathered
\La=\R[\bp_1,\bp_2,\bp_3\dots]\, \to\,
\Alg=\R[\tp_2,\tp_3,\dots],\\
\bp_1\to0, \qquad \bp_k\to\tp_k, \quad k=2,3,\dots
\endgathered \tag2.18
$$
This definition is motivated by Proposition 2.11 below.
\endexample

This morphism induces an algebra isomorphism $\La/\bp_1\La\to\Alg$.
Let us emphasize the difference from the isomorphism $\La\to\Alg$
(Definition 1.6).

The weight grading induces a filtration in $\Alg$, which we call the
{\it weight filtration\/} and denote by the same symbol
$\wt(\,\cdot\,)$. Note that
$$
\wt(p_k)=k+1, \qquad k=1,2,\dots,
$$
because the top weight homogeneous component of $p_k$ is
$\tp_{k+1}/(k+1)$, see Corollary 2.8.

\example{Definition 2.10} a) We define an action of the multiplicative
group of positive real numbers on the set $\Dz$  by setting
$$
\om^s(x)=s^{-1}\om(sx),
\qquad \om\in\Dz, \quad s>0, \quad x\in\R.
$$
In other words, the graph of $y=\om^s(x)$ is obtained from that of
$y=\om(x)$ by the transformation $(x,y)\mapsto(s^{-1}x,s^{-1}y)$.

b) Since $\Alg=\R[\tp_2,\tp_3,\dots]$, we may define the symbol
$f[\om]$ (where $\om$ ranges over $\Dz$) for any $f\in\Alg$.
Specifically, write $f$ as a polynomial in $\tp_2,\tp_3,\dots$ and
then specialize each $\tp_k$ to $\tp_k[\om]$. In this way, we realize
$\Alg$ as an algebra of functions on $\Dz$.
\endexample

\proclaim{Proposition 2.11} Let $f\in\Alg$ be homogeneous with
respect to the weight grading, Definition 2.9. Then for any
$\om\in\Dz$ and $s>0$,
$$
f[\om^s]=s^{-\wt(f)}f[\om].
$$
\endproclaim

\demo{Proof} By the definition of the weight grading, it suffices to
check that
$$
\tp_k[\om^s]=s^{-k}\tp_k[\om],
\qquad k=2,3,\dots, \quad \om\in\Dz, \quad s>0. \tag2.19
$$
Remark that the function $\si(x)=\tfrac12(\om(x)-|x|)$ transforms in
the same way as $\om(x)$. Then \tht{2.19} is clear from Proposition 2.2.\qed
\enddemo

\head \S3. The elements $\pp_k$ \endhead

Let $\frak S_n$ be the symmetric group of degree $n$.  Recall that
both irreducible characters and conjugacy classes of $\frak S_n$ are
indexed by the same set, the set of partitions of $n$ or,
equivalently, of Young diagrams with $n$ boxes. We denote this set by
$\Y_n$\,. For $\la, \rho\in\Y_n\,$, we denote
by $\chi^\la$ the irreducible character of $\frak S_n$ indexed by
$\la$, and by $\chi^\la_\rho$ the value of $\chi^\la$ on the
conjugacy class indexed by $\rho$.

In particular, the partition $\rho=(1^n)=(1,\dots,1)$ corresponds to
the trivial conjugacy class $\{e\}\subset\frak S_n$\,, so that
$\chi^\la_{(1^n)}$ equals the dimension of $\chi^\la$; we denote this
number by $\dim\la$.

\example{Definition 3.1} For $k=1,2,\dots$, let $\pp_k$ be following
function on $\Y$:
$$
\pp_k(\la)=\cases n^{\fd k}\,
\cdot\,\dfrac{\chi^\la_{(k,1^{n-k})}}{\dim\la}\,,
& n:=|\la|\ge k,\\
0, & n<k,
\endcases
$$
where
$$
n^{\fd k}=n(n-1)\dots(n-k+1)
$$
and
$$
(k,1^{n-k})=(k, 1,\dots,1)\in\Y_n\,.
$$
\endexample

\proclaim{Proposition 3.2} For any $k=1,2,\dots$ and any $\la\in\Y$,
$\pp_k(\la)$ equals the coefficient of $z^{-1}$ in the expansion of
the function
$$
-\frac1k\,(z-\tfrac12)^{\fd k}\,\frac{\Phi(z;\la)}{\Phi(z-k;\la)} \tag3.1
$$
in descending powers of $z$ about the point $z=\infty$.
\endproclaim

\demo{Proof} First, assume $n<k$. Then, by the definition,
$\pp_k(\la)=0$, and we have to prove that the coefficient in question
equals 0, too. It suffices to prove that \tht{3.1} is a polynomial in
$z$. By Proposition 1.2, \tht{3.1} equals
$$
-\frac1k\,(z-\tfrac12)^{\fd k}\,
\prod_{i=1}^d\frac{z+b_i}{z-a_i}\,\cdot\,
\prod_{j=1}^d\frac{z-a_j-k}{z+b_j-k}\,.
$$
Note that $a_i\ne k-b_j$ for any $i,j=1,\dots,d$, because
$a_i+b_j\le n<k$. Therefore, all the factors $z-a_i$ and
$z+b_j-k=z-(k-b_j)$ are pairwise distinct. Each of them cancels with
one of the factors in the product $(z-\tfrac12)^{\fd k}$, because
$$
a_i, \, k-b_j\,\in\, \{\tfrac12,\tfrac32,\dots,k-\tfrac12\}.
$$
This concludes the proof in the case $n<k$.

Now we shall assume $n\ge k$. Then we use a formula due to Frobenius
(see \cite{Ma, Ex. I.7.7}) which says that $\pp_k(\la)$ equals the
coefficient of $z^{-1}$ in the expansion of the function
$$
F(z)=-\frac1k\, z^{\fd k}\,
\prod_{i=1}^n \frac{z-\la_i-n+i-k}{z-\la_i-n+i}
$$
about $z=\infty$. In other words,
$$
\pp_k(\la)=-\underset{{z=\infty}}\to{\Res}(F(z)dz).
$$
After simple transformations we get
$$
F(z)=-\frac1k\,(z-n)^{\fd k}\,
\frac{\Phi(z-n+\tfrac12;\la)}{\Phi(z-n+\tfrac12-k;\la)}\,.
$$
The residue at $z=\infty$ will not change under the shift $z\mapsto
z+n-\tfrac12$. Consequently,
$$
\pp_k(\la)=-\underset{{z=\infty}}\to{\Res}\left(
-\frac1k\,(z-\tfrac12)^{\fd k}\,
\frac{\Phi(z;\la)}{\Phi(z-k;\la)}\right),
$$
which completes the proof. \qed
\enddemo

We shall employ the following notation. Given a formal series $A(t)$,
let
$$
[t^k]\{A(t)\}=\, \langle\text{the coefficient of $t^k$ in $A(t)$}\rangle.
$$
The next result is due to Wassermann \cite{Wa, \S III.6}.

\proclaim{Proposition 3.3} For any $k=1,2,\dots$, the function
$\pp_k(\la)$ introduced in Definition 2.1 belongs to the algebra
$\Alg$. Its expression through the generators $p_1,p_2,\dots$ of
$\Alg$ can be described as follows:
$$\pp_k=[t^{k+1}]\left\{
-\frac1k\, \prod_{j=1}^k(1-(j-\tfrac12)t)\,\cdot\,
\exp\left(\sum_{j=1}^\infty\frac{p_jt^j}j\,(1-(1-kt)^{-j})\right)\right\}.
\tag3.2
$$
\endproclaim

\demo{Proof}
By Proposition 3.2,
$$
\pp_k(\la)=[t^{k+1}]\left\{
-\tfrac1k\,t^k\,(t^{-1}-\tfrac12)^{\fd k}\,
\dfrac{\Phi(t^{-1};\la)}
{\Phi\left(\left(\frac t{1-kt}\right)^{-1};\la\right)}\right\}.
$$
We have
$$
t^k(t^{-1}-\tfrac12)^{\fd k}=\prod_{j=1}^k(1-(j-\tfrac12)t)
$$
and, by Proposition 1.4,
$$
\Phi(t^{-1};\la)=\exp\left(\sum_{j=1}^\infty\frac{p_j(\la)}j\,t^j\right).
$$
This yields \tht{3.2}, which in turn implies that $\pp_k\in\Alg$.\qed
\enddemo

The expression \tht{3.2} can be written in the form
$$
\gather
\pp_k=-\frac1k\,[t^{k+1}]\left\{
(1+\ep_0(t))\exp\left(-\sum_{j=1}^\infty
kp_jt^{j+1}(1+\ep_j(t))\right)\right\}\\
=-\frac1k\,[t^{k+1}]\left\{(1+\ep_0(t))\sum_{m=0}^\infty
\frac{(-1)^m}{m!}\left(\sum_{j=1}^\infty
kp_jt^{j+1}(1+\ep_j(t))\right)^m\right\}. \tag3.3
\endgather
$$
Here each $\ep_r(t)$ is a power series of the form
$c_1t+c_2t^2+\dots$, where the coefficients $c_1,c_2,\dots$ do not
involve the generators $p_1,p_2,\dots$.

Using \tht{3.3} we can readily evaluate the top homogeneous component of
$\pp_k$ with respect both to the canonical grading and the weight
grading in $\Alg$.

\proclaim{Proposition 3.4} In the canonical grading, see Definition
1.6, the highest term of $\pp_k$ equals $p_k$.
\endproclaim

\demo{Proof} Apply \tht{3.3} and write the expression in the curly
brackets as a sum of terms of the form $\const\cdot P\cdot t^r$,
where $P$ stands for a monomial in $p_1,p_2,\dots$. We search for
terms with $r=k+1$ and such that $\deg P$, the total degree of $P$,
counted with the convention that $\deg p_k=k$, is maximal possible.

The first observation is that all terms $\const\cdot P\cdot t^r$
involving at least one factor coming from $\ep_0(t), \ep_1(t),\dots$
are negligible, because the epsilon factors diminish the difference
$\deg P-r$. Removing $\ep_0(t), \ep_1(t),\dots$, we get
$$
\pp_k=-\frac1k\,[t^{k+1}]
\left\{1+\sum_{m=1}^\infty\frac{(-1)^m}{m!}
\left(\sum_{j=1}^\infty kp_jt^{j+1}\right)^m\right\}+\dots,
$$
where dots stand for lower degree terms. The summand with $m=1$ has a
unique term with $r=k+1$. This term this $-kp_kt^{k+1}$, and its
contribution is $p_k$.

The second observation is that the summands with $m=2,3,\dots$ are
negligible, because, in the corresponding terms, $\deg P-r=-m\le-2$,
so that $r=k+1$ implies $\deg P<k$.

We conclude that $\pp_k=p_k+\dots$\,. \qed
\enddemo

\proclaim{Proposition 3.5} Let $k=1,2,\dots$\,. In the weight
grading, the top homogeneous component of $\pp_k$ has weight $k+1$
and can be written as
$$
\gather
-\frac1k\,[t^{k+1}]\left\{\exp\big(
-k\sum_{j=2}^\infty \frac{\tp_j}j\,t^j\big)\right\}\tag3.4\\
=\frac{\tp_{k+1}}{k+1}\,+\,\langle\text{\rm a homogeneous polynomial in
$\tp_2,\dots,\tp_k$ of total weight $k+1$}\rangle. \tag3.5
\endgather
$$
\endproclaim

\demo{Proof} Apply \tht{3.3} and recall that $\wt(p_j)=j+1$, because
$$
p_j=\frac{\tp_{j+1}}{j+1}\,+\,\langle\text{a linear combination of
$\tp_2,\dots,\tp_j$}\rangle.
$$
As in the situation of Proposition 3.4, we may neglect the epsilon
factors, which affect only terms of lower weight. For the same reason,
we may replace each $p_j$ by $\tp_{j+1}/(j+1)$. This leads to \tht{3.4},
and \tht{3.5} follows from \tht{3.4}. \qed
\enddemo

In Proposition 3.7 we invert the result of Proposition 3.5. Beforehand we
state the following general fact.

\proclaim{Proposition 3.6} Let $a_2,a_3,\dots$ and $b_2,b_3,\dots$ be
two families of elements in a commutative algebra. Let
$$
A(t)=1+\sum_{j=2}^\infty a_jt^j, \qquad
B(u)=1+\sum_{j=2}^\infty b_ju^j
$$
be their generating series, and set
$$
\tA(t)=\ln A(t)=\sum_{j=2}^\infty\frac{\ta_j}{j}\,t^j\,.
$$

Then the following conditions are equivalent:
\medskip

{\rm (i)} The formal transformations $x\to xA(x)$ and $x\to x/B(x)$
are inverse to each other.
\medskip

{\rm (ii)} $b_k=-\,\frac1{k-1}\,[t^k]\{A^{-(k-1)}(t)\}, \quad
k=2,3,\dots\,$.
\medskip

{\rm(iii)} $a_k=\frac1{k+1}\,[u^k]\{B^{k+1}(u)\}, \quad
k=2,3,\dots\,$.
\medskip
{\rm(iv)} $\ta_k=[u^k]\{B^k(u)\}, \quad k=2,3,\dots\,$.
\endproclaim

\demo{Proof} This is a variation of Lagrange's inversion formula and
can be proved by the standard argument, see, e.g., \cite{Wi},
\cite{Ma, Example I.2.24}. \qed
\enddemo

In Proposition 3.7 we use only a part of the claims of Proposition
3.6. Another part will be used later on.

\proclaim{Proposition 3.7} For $k=2,3,\dots$
$$
\gather
\tp_k=[u^k]\left\{\big(1+\sum_{j=2}^\infty
\pp_{j-1}u^j\big)^k\right\}\,+\,\dots\\
=\sum\Sb m_2,m_3,\dots\\ 2m_2+3m_3+\dots=k\endSb
\frac{k^{\fd\sum m_i}}{\prod m_i!}\,
\prod_{i\ge2}(\pp_{i-1})^{m_i}\,+\,\dots, \tag3.6
\endgather
$$
where dots mean a polynomial in $\pp_1,\pp_2,\dots,\pp_{k-2}$ of total
weight $\le k-1$, where $\wt(\pp_i)=i+1$.
\endproclaim

\demo{Proof} Assume that $\ta_2,\ta_3,\dots$ and $b_2,b_3,\dots$ are
elements of a commutative algebra such that
$$
b_k=\frac1k\,\ta_k+X_k(\ta_2,\dots,\ta_{k-1}), \quad k=2,3,\dots,\tag3.7
$$
where $X_k$ is an inhomogeneous polynomial in $k-2$ variables
such that $\wt(X_k)\le k$, where $\wt(X_k)$ denotes the total weight
counted with the convention that $\wt(\ta_j)=j$.

Then, as is readily seen,
$$
\frac1k\,\ta_k=b_k+Y_k(b_2,\dots,b_{k-1}), \quad k=2,3,\dots,\tag3.8
$$
where, likewise, $Y_k$ is a polynomial of total weight $\le k$, with
the convention that $\wt(b_j)=j$.

Moreover, the top weight homogeneous component of $Y_k$ depends only
on the top weight homogeneous components of $X_2,\dots,X_k$.

Now let us set
$$
\gather
\ta_k=\tp_k\,, \quad
b_k=\pp_{k-1}\,, \qquad k=2,3,\dots,\\
\tA(t)=\sum_{j=2}^\infty \frac{\ta_j}j t^j, \quad
A(t)=\exp\tA(t), \quad
B(u)=1+\sum_{j=2}^\infty b_j u^j.
\endgather
$$
By Proposition 3.5, we have
$$
b_k=-\,\frac1{k-1}\,[t^k]\{A^{-(k-1)}(t)\}+\dots, \quad k=2,3,\dots,
$$
where dots mean terms of lower weight. By \tht{3.6}, these relations are
of the form \tht{3.7}. Therefore, to evaluate the inverse relations
\tht{3.8} up to lower weight terms, we may use formula \tht{iv} of
Proposition 3.6. This yields
$$
\ta_k=[u^k]\{B^k(u)\}+\dots,\quad k=2,3,\dots,
$$
which is exactly \tht{3.6}. \qed
\enddemo

\head \S4. The basis $\{\pp_\rho\}$ and filtrations in $\Alg$
\endhead

In this section we review some results of \cite{KO} and \cite{IK}.

\example{Definition 4.1} To any partition $\rho$ we assign a function
$\pp_\rho$ on $\Y$ as follows. Let $r=|\rho|$, let $\la\in\Y$, and
denote $n=|\la|$. Then
$$
\pp_\rho(\la)=\cases n^{\fd r}\cdot \dfrac{\chi^\la_{\rho\cup1^{n-r}}}
{\dim\la}\,, & n\ge r,\\
0, & n<r,
\endcases
$$
where $\rho\cup1^{n-r}=(\rho,1,\dots,1)\in\Y_n$.
\endexample

When $\rho$ consists of a single part, $\rho=(r)$, then this reduces
to Definition 3.1.

Given a partition $\rho$, we shall denote by $m_i=m_i(\rho)$ the
multiplicity of $i$ in $\rho$:
$$
m_i(\rho)=\operatorname{Card}\{j\mid \rho_j=i\}, \quad i=1,2,\dots\,.
$$
By $\ell(\rho)$ we denote the number of nonzero parts of $\rho$ (the
length of $\rho$). We have $\ell(\rho)=\sum m_i(\rho)$.

Similarly to the conventional notation for the algebra $\La$, we set
$$
p_\rho=p_{\rho_1}\dots p_{\rho_{\ell(\rho)}}=\prod_{i}
p_i^{m_i(\rho)}\,.
$$
The elements $p_\rho$ form a homogeneous (in the canonical grading)
basis in the algebra $\Alg$. Note that $\deg p_\rho=|\rho|$.

The next result generalizes Proposition 3.3 (first claim) and
Proposition 3.4. It was first announced in \cite{VeK2}.

\proclaim{Proposition 4.2} For any partition $\rho$, the function
$\pp_\rho$ introduced in Definition 4.1 is an element of $\Alg$. In
the canonical grading, the top degree homogeneous component of
$\pp_\rho$ equals $p_\rho$.
\endproclaim

\demo{Proof} Different proofs are given in \cite{KO} and \cite{OO}.
See also \cite{IK}, \cite{LaT}, \cite{ORV}. \qed
\enddemo

\proclaim{Corollary 4.3} The elements $\pp_\rho$ form a basis in
$\Alg$.
\endproclaim

Note that this basis is inhomogeneous both in the canonical grading
and the weight grading.

Given two partitions $\si,\tau$, we denote by $\si\cup\tau$ the
partition obtained by joining the parts of both partitions and then
arranging them in descending order. In other words, $\si\cup\tau$ is
characterized by
$$
m_i(\si\cup\tau)=m_i(\si)+m_i(\tau), \quad i=1,2,\dots\,.
$$

\proclaim{Corollary 4.4} For any partitions $\si,\tau$,
$$
\pp_\si \pp_\tau=\pp_{\si\cup\tau}+\dots,
$$
where dots mean lower degree terms with respect to the canonical
grading.
\endproclaim

Here and in what follows we define the degree of an inhomogeneous
element as the maximal degree of its nonzero homogeneous components.
In other words, we switch from the grading to the corresponding
filtration.

Later on it will be shown that the claim of Corollary 4.4 also holds for
the weight grading (or filtration), see Propositions 4.9 and 4.10.

Let $f^\rho_{\si\tau}$ denote the structure constants of the algebra
$\Alg$ in the basis $\{\pp_\rho\}$. I.e.,
$$
\pp_\si\pp_\tau=\sum_{\rho} f^\rho_{\si\tau}\, \pp_\rho\,.
$$
By Corollary 4.4, $f^\rho_{\si\tau}\ne0$ implies
$|\rho|\le|\si|+|\tau|$. Moreover,
$$
f^{\si\cup\tau}_{\si\tau}=1. \tag4.1
$$

Recall the conventional notation \cite{Ma, \S I.2}
$$
z_\rho=\prod_i i^{m_i(\rho)}\, m_i(\rho)!
$$

\proclaim{Proposition 4.5} Let $\rho,\si,\tau$ be arbitrary
partitions. We have
$$
f^\rho_{\si\tau}=\frac{z_\si z_\tau}{z_\rho}\, g^\rho_{\si\tau}\,,
$$
where $g^\rho_{\si\tau}$ can be evaluated as follows.

Fix a set $X$ of cardinality $|\rho|$ and a permutation $s: X\to X$
whose cycle structure is given by $\rho$. Then $g^\rho_{\si\tau}$
equals the number of quadruples $(X_1,s_1,X_2,s_2)$ such that:
\medskip

{\rm(i)} $X_1\subseteq X, \quad X_2\subseteq X, \quad X_1\cup X_2=X$.
\medskip

{\rm(ii)} $|X_1|=|\si|$ and $s_1: X_1\to X_1$ is a permutation of cycle
structure $\si$.
\medskip

{\rm(iii)} Likewise, $|X_2|=|\tau|$ and $s_2: X_2\to X_2$ is a
permutation of cycle structure $\tau$.
\medskip

{\rm(iv)} Denote by $\bar s_1: X\to X$ and $\bar s_2: X:\to X$ the
natural extensions of $s_{1,2} $ from $X_{1,2}$ to the whole $X$.
I.e., $\bar s_{1,2}$ is trivial on $X\setminus X_{1,2}$. Then the
condition is that $\bar s_1 \bar s_2=s$.
\endproclaim

\demo{Proof} See \cite{IK, Proposition 6.2 and Theorem 9.1}. \qed
\enddemo

\example{Definition 4.6} Fix an arbitrary subset $J\subseteq\N$,
where $\N=\{1,2,\dots\}$. For any partition $\rho$, set
$$
|\rho|_J=|\rho|+\sum_{j\in J} m_j(\rho).
$$
In particular, in the two extreme cases we have
$$
|\rho|_\varnothing=|\rho|, \quad |\rho|_\N=|\rho|+\ell(\rho).
$$
Next, following \cite{IK}, define a filtration of the vector space
$\Alg$ by setting
$$
\deg_J(\pp_\rho)=|\rho|_J
$$
and, more generally, for any $f=\sum_\rho f_\rho \pp_\rho\in\Alg$,
$$
\deg_J(f)=\max_{\rho:\, f_\rho\ne0} |\rho|_J\,.
$$
\endexample

\proclaim{Proposition 4.7} For any $J\subseteq\N$, the filtration by
$\deg_J(\,\cdot\,)$ as defined above is compatible with the
multiplication in $\Alg$. I.e., for any partitions $\rho,\si,\tau$,
$$
f^\rho_{\si\tau}\ne0 \, \Longrightarrow \,
|\rho|_J\le|\si|_J+|\tau|_J\,,
$$
so that this is an algebra filtration.
\endproclaim

\demo{Proof} The argument presented below is a slightly rewritten
version of that given in \cite{IK, Proposition 10.3}.

Assume we are given partitions $\rho,\si,\tau$ such that
$f^\rho_{\si\tau}\ne0$. Fix a set $X$ and a permutation $s:X\to X$ as
in the statement of Proposition 4.5. By that proposition, there
exists a quadruple $\{X_1,s_1,X_2,s_2\}$ satisfying the four
conditions (i)--(iv). Fix any such quadruple.

Decompose each of the permutations $s,s_1,s_2$ into cycles and denote
by $C_J(\,\cdot\,)$ the set of all cycles whose lengths belong to the
set $J$. Write
$$
C_J(s_1)=A_J(s_1)\sqcup B_J(s_1), \quad
C_J(s_2)=A_J(s_2)\sqcup B_J(s_2), \tag4.2
$$
where $A_J(s_1)\subseteq C_J(s_1)$ denotes the subset of those cycles
of $s_1$ that are entirely contained in $X_1\setminus X_2$, while
$B_J(s_1)\subseteq C_J(s_1)$ denotes the subset of those cycles of
$s_1$ that have a nonempty intersection with $X_1\cap X_2$. (Note that
we count fixed points viewed as cycles of length 1, provided
that $1\in J$.)
The sets $A_J(s_2)$ and $B_J(s_2)$ are defined similarly.

In this notation we have
$$
C_J(s)=A_J(s_1)\sqcup A_J(s_2)\sqcup B_J(s), \tag4.3
$$
where $B_J(s)$ is the set of those cycles in $C_J(s)$ that intersect both
$X_1$ and $X_2$.

Remark that
$$
|\rho|_J=|X|+|C_J(s)|, \quad
|\si|_J=|X_1|+|C_J(s_1)|, \quad
|\tau|_J=|X_2|+|C_J(s_2)|.
$$
Therefore, the required inequality $|\rho|_J\le |\si|_J+|\tau|_J$
means
$$
|X|+|C_J(s)|\le|X_1|+|C_J(s_1)|+|X_2|+|C_J(s_2)|.
$$
By virtue of \tht{4.2} and \tht{4.3} this is equivalent to
$$
|X|+|B_J(s)|\le |X_1|+|B_J(s_1)|+|X_2|+|B_J(s_2)|. \tag4.4
$$
We shall establish a stronger inequality,
$$
|X|+|B_J(s)|\le |X_1|+|X_2|, \tag4.5
$$
which is equivalent to
$$
|B_J(s)|\le |X_1\cap X_2|. \tag4.6
$$
To prove the latter inequality, we shall show that each cycle
$c\in B_J(s)$ contains a point of $X_1\cap X_2$.

By the definition of $B_J(s)$, $c$ contains both points of $X_1$ and
of $X_2$. Therefore, there exist points $x_1\in X_1\cap c$ and
$x_2\in X_2\cap c$ such that $sx_1=x_2$. We claim that either $x_1$ or
$x_2$ lies in $X_1\cap X_2$.  Indeed, if $x_1\in X_1\setminus X_2$
then
$$
x_2=sx_1=\bar s_1 \bar s_2 x_1=\bar s_1 x_1=s_1 x_1\in X_1 \,.
$$
This shows that $x_2\in X_1\cap X_2$, which completes the proof. \qed
\enddemo

\proclaim{Corollary 4.8 (of the proof)} Let $f^\rho_{\si\tau}\ne0$
and $|\rho|_J=|\si|_J+|\tau|_J$. Then, in the notation of the proof of
Proposition 4.7, $B_J(s_1)=\varnothing$, $B_J(s_2)=\varnothing$, and
\tht{4.6} is actually an equality.
\endproclaim

\demo{Proof} Indeed, the equality $|\rho|_J=|\si|_J+|\tau|_J$ means
that \tht{4.4} is an equality. Then \tht{4.5} is an equality, too.
This implies all the claims. \qed
\enddemo

\proclaim{Proposition 4.9} Assume $J=\N$. For any partitions
$\si,\tau$,
$$
\pp_\si\pp_\tau=\pp_{\si\cup\tau}\,+\,
\langle\text{\rm a linear combination of
$\pp_\rho$'s with $|\rho|_\N<|\si|_\N+|\tau|_\N$}\rangle.
$$
\endproclaim

\demo{Proof} We have
$$
\pp_\si\pp_\tau=\sum_\rho f^\rho_{\si\tau}\,\pp_\rho\,.
$$
By Proposition 4.7, only partitions $\rho$ with
$|\rho|_\N\le|\si|_\N+|\tau|_\N$ can really contribute.

By Corollary 4.8, if $f^\rho_{\si\tau}\ne0$ and
$|\rho|_\N=|\si|_\N+|\tau|_\N$\,, then both $B_\N(s_1)$ and $B_\N(s_2)$
are empty, which implies $X_1\cap X_2=\varnothing$. Therefore,
$\rho=\si\cup\tau$. Finally, by formula \tht{4.1}, which we have derived
from Corollary 4.4, $f^{\si\cup\tau}_{\si\tau}=1$.
This completes the proof. \qed
\enddemo

Note that formula \tht{4.1} can also be obtained from Proposition 4.9.

\proclaim{Proposition 4.10} The filtration of $\Alg$ defined by
$\deg_\N(\,\cdot\,)$ coincides with the weight filtration.
\endproclaim

\demo{Proof} For any $r=0,1,\dots$, let $\Alg'_r\subset\Alg$ denote
the $r$th member of the first filtration, and let
$\Alg''_r\subset\Alg$ has the same meaning for the second
filtration. Recall that
$$
\gather
\Alg'_r=\Span\{\pp_\rho\mid |\rho|_\N=|\rho|+\ell(\rho)\le r\},\\
\Alg''_r=\Span\{\tp_{k_1}\dots \tp_{k_l}\mid
k_1,\dots,k_l\ge2, \quad k_1+\dots +k_l\le r\}.
\endgather
$$

Clearly,
$$
\Alg'_0=\Alg''_0=\Alg'_1=\Alg''_1=\R\cdot1.
$$
We shall prove that for any $r\ge2$, both $\Alg'_r\subseteq \Alg''_r$
and $\Alg''_r\subseteq \Alg'_r$. By Proposition 4.9,
$$
\Alg'_r=\Span\{\pp_{k_1}\dots\pp_{k_l}\mid
k_1,\dots k_l\ge 1, \quad k_1+\dots+k_l+l\le r\}.
$$
Therefore, it suffices to show that
$$
\gather
\pp_k\in\Span\{\tp_{k_1}\dots\tp_{k_l}\mid
k_1,\dots,k_k\ge2, \quad k_1+\dots+k_l\le k+1\}, \\
\tp_k\in\Span\{\pp_{k_1}\dots\pp_{k_l}\mid
k_1,\dots, k_l\ge1, \quad k_1+\dots+k_l+l\le k\}.
\endgather
$$
The first inclusion follows from Proposition 3.5, and the second
inclusion follows from Proposition 3.7. \qed
\enddemo

In the remaining part of the section we focus on the filtration
corresponding to $J=\{1\}$. It first appeared in \cite{Ke1}, and we
propose to call it the {\it Kerov filtration\/} of $\Alg$.

Let us abbreviate
$$
|\rho|_1=|\rho|_{\{1\}}=|\rho|+m_1(\rho), \quad
\deg_1(\,\cdot\,)=\deg_{\{1\}}(\,\cdot\,).
$$
The next three results will be used in \S6.

\proclaim{Proposition 4.11} For any partition $\si$,
$$
\pp_\si\pp_1=\pp_{\si\cup1}\,+\,
\langle\text{\rm a term of lower degree with
respect to $\deg_1(\,\cdot\,)$}\rangle.
$$
\endproclaim

\demo{Proof} Actually, the following exact formula holds:
$$
\pp_\si\pp_1=\pp_{\si\cup1}+|\si|\cdot\pp_\si\,. \tag4.7
$$
To prove this, apply Definition 4.1 and evaluate both sides at
a partition $\la$. It suffices to assume that $n=|\la|$ is large
enough, $n\ge|\si|+1$. Then we get, abbreviating $k=|\si|$,
$$
\pp_\si(\la)
=n^{\fd k}\cdot\frac{\chi^\la_{\si\cup1^{n-k}}}{\dim\la}\,,
\quad
\pp_{\si\cup1}(\la)
=n^{\fd(k+1)}\cdot\frac{\chi^\la_{\si\cup1^{n-k}}}{\dim\la}\,,
$$
because $(\si\cup1)\cup1^{n-k-1}=\si\cup1^{n-k}$.

Therefore, the verification of \tht{4.7} at $\la$ reduces to the
relation
$$
n^{\fd k}\cdot n=n^{\fd(k+1)}+n^{\fd k}\cdot k,
$$
which is trivial. \qed
\enddemo

Note that Proposition 4.11 can also be obtained from Proposition 4.7.
We shall use this approach in the next proposition.

\proclaim{Proposition 4.12} For any partition $\si$ and any $k\ge2$,
$$
\pp_\si\pp_k=\pp_{\si\cup k}+
\cases k\cdot m_k(\si)\cdot \pp_{(\si\setminus k)\cup1^k}\,
+\,\dots, & m_k(\si)\ge1, \\
\dots, & m_k(\si)=0,
\endcases \tag4.8
$$
where the partition $\si\setminus k$ is obtained from $\si$ by
removing one part equal to $k$, i.e.,
$$
m_i(\si\setminus k)=\cases m_i(\si), & i\ne k, \\
m_k(\si)-1, & i=k,
\endcases
$$
and dots mean terms of lower degree, i.e., with
$\deg_1(\,\cdot\,)<|\si|_1+k$.
\endproclaim

\demo{Proof} Assume first that $\tau$ is an arbitrary partition (not
necessarily $\tau=(k)$) and search for partitions $\rho$ such that
$f^\rho_{\si\tau}\ne0$ and $|\rho|_1=|\si|_1+|\tau|_1$.  Let us
employ the notation introduced in the proof of Proposition 4.7 and
apply Corollary 4.8. We get $B_{\{1\}}(s_1)=\varnothing$,
$B_{\{1\}}(s_2)=\varnothing$, and $|B_{\{1\}}(s)|=|X_1\cap X_2|$. This
means that in $X_1\cap X_2$ there is no 1-cycle ($=$fixed point) for $s_1$
and $s_2$, but all points of $X_1\cap X_2$ are are fixed by $s$.

This shows that either $X_1\cap X_2=\varnothing$ or $X_1\cap X_2$
entirely consists of common nontrivial cycles of the permutations
$s_1$ and $s_2^{-1}$.

Now apply this conclusion to the special case $\tau=(k)$ that we
need. Recall that $k\ge2$. The first possibility, $X_1\cap
X_2=\varnothing$, means that $\rho=\si\cup\tau=\si\cup(k)$. Then the
corresponding coefficient $f^\rho_{\si\tau}$ equals 1, see \tht{4.1}.
This explains the term $\pp_{\si\cup k}$ in the right--hand side of
\tht{4.8}.

The second possibility means that $X_1\supseteq X_2$\,, because
$s_2^{-1}$ reduces to a single $k$--cycle, which is also a $k$--cycle
of $s_1$.  This implies that $m_k(\si)\ge1$ and
$\rho=(\si\setminus k)\cup1^k$.

It remains to evaluate the coefficient $f^\rho_{\si\tau}$. Let us
abbreviate $m=m_k(\si)$, $l=m_1(\si)$. We must prove that
$f^\rho_{\si\tau}=km$. To do this we apply Proposition 4.5. We
readily get
$$
\frac{z_\rho}{z_\si z_\tau}=\frac{(k+l)!}{l!k^2m}\,,
$$
so that $f^\rho_{\si\tau}=km$ is equivalent to
$$
g^\rho_{\si\tau}=\frac{(k+l)!}{l!k}\,.
$$

Let us check this. By the definition of $g^\rho_{\si\tau}$\,, in our
situation it equals the number of ways to choose a $k$--cycle inside
a $(k+l)$--point set. This number equals
$$
\frac{(k+l)!}{k!l!}\cdot(k-1)!=\frac{(k+l)!}{l!k}\,,
$$
(the number of $k$--point subsets inside a $(k+l)$--point set, times
the number of different $k$--cycle structures on a given $k$--point
set). This concludes the proof. \qed
\enddemo

Finally, note that
$$
\inv(\pp_\rho)=(-1)^{|\rho|+\ell(\rho)}\pp_\rho\,. \tag4.9
$$
Indeed, this follows from the definition of $\pp_\rho$. In particular,
$$
\inv(\pp_k)=(-1)^{k-1}\pp_k\,. \tag4.10
$$
This symmetry property will be used in the proofs of Proposition 7.3
and Theorem 10.2.

\proclaim{Corollary 4.13 (of the proof)} Let $\si$ and $\tau$ be two
partitions with no common part, i.e., for any $i=1,2,\dots$, at least
one of the multiplicities $m_i(\si)$, $m_i(\tau)$ vanishes. Then
$$
\pp_\si \pp_\tau=\pp_{\si\cup\tau}\,+\,
\langle\text{\rm terms of lower degree $\deg_1(\,\cdot\,)$}\rangle.
$$
\endproclaim

\demo{Proof} Let $\rho$, $X_1$, $X_2$ be as in the beginning of the
proof of Proposition 4.12. Recall the claim stated in the second
paragraph of that proof: either $X_1\cap X_2=\varnothing$ or
$X_1\cap X_2$ entirely consists of common nontrivial cycles of the
permutations $s_1$ and $s_2^{-1}$. The second possibility contradicts
the assumption that $\si$ and $\tau$ have no common part. Hence the
first possibility holds, which means that $\rho=\si\cup\tau$. We know
that the corresponding coefficient $f^\rho_{\si\tau}$ equals 1, which
concludes the proof. \qed
\enddemo

\head \S5. The Plancherel measure and the law of large numbers
\endhead

Consider the set $\Y_n$ of Young diagrams with $n$ boxes,
$n=1,2,\dots$, and equip it with the measure $M_n$, defined by
$$
M_n(\la)=\frac{\dim^2\la}{n!}\,, \quad \la\in\Y_n\,.
$$
This is a probability measure, because, by Burnside's theorem,
$$
\sum_{\la\in\Y_n}\dim^2\la=|\frak S_n|=n!.
$$
It is called the {\it Plancherel measure,\/} see \cite{VeK1}, \cite{VeK2},
\cite{VeK3} for more details.

Given a function $f$ on $\Y_n$\,, we define by $\m f$ its expectation
with respect to $M_n$. That is,
$$
\m f=\sum_{\la\in\Y_n}f(\la)M_n(\la)
=\frac1{n!}\,\sum_{\la\in\Y_n}f(\la)\dim^2\la.
$$

If $f$ is a function on the whole set $\Y$, we write $\m f$ instead
of $\m {f\mid_{\Y_n}}$\,. We shall use this convention for functions
$f\in\Alg$. In this way we get the family of linear
functionals $\m{\,\cdot\,}$\,, $n=1,2,\dots$, on the algebra $\Alg$.
These functionals have a very simple form on the basis
$\{\pp_\rho\}$.

\proclaim{Proposition 5.1} For any partition $\rho$,
$$
\m{\pp_\rho}=\cases n^{\fd r}, & \rho=(1^r), \quad r=1,2,\dots, \\
0, & \rho\ne(1^r).
\endcases
$$
\endproclaim

\demo{Proof} Set $r=|\rho|$. If $n<r$ then $\pp_\rho$ vanishes on
$\Y_n$, which agrees with the formula in question, because $n^{\fd
r}=0$ whenever $n<r$.

Assume $n\ge r$. By the definition of $\pp_\rho$\,, see Definition 4.1,
$$
\m{\pp_\rho}=n^{\fd r}\cdot \frac1{n!}\cdot
\sum_{\la\in\Y_n}\chi^\la_{\rho\cup1^{n-r}}\,\dim\la.
$$
Remark that the sum above equals the value of the regular character
(i.e., the character of the regular representation of $\frak S_n$) on
the conjugacy class $\rho\cup1^{n-r}$. But the regular character is
the delta function at $\{e\}\subset\frak S_n$\,, multiplied by $n!$.
Therefore, the sum in question vanishes unless $\rho\cup1^{n-r}$ is
the trivial class (i.e., $\rho$ itself is trivial, $\rho=(1^r)$), in
which case the sum equals $n!$. Consequently, we get $n^{\fd r}$. \qed
\enddemo

\proclaim{Proposition 5.2} For any $f\in\Alg$, the expectation $\m f$
is a polynomial in $n$. The degree of this polynomial is bounded from
above by $\tfrac12\deg_1(f)$.
\endproclaim

\demo{Proof} It suffices to check this for $f=\pp_\rho$. Let
$r=|\rho|$. If $\rho\ne(1^r)$ then $\m f\equiv0$, by virtue of
Proposition 5.1. If $\rho=(1^r)$ then, by Proposition 5.1, $\m
f=n^{\fd r}$. On the other hand, $\deg_1(f)=2r$, which agrees with
the claim. \qed
\enddemo

We define the function $\Om(x)$ on $\R$ by
$$
\Om(x)=\cases \frac2\pi(x\arcsin\tfrac x2+\sqrt{4-x^2}), & |x|\le2,\\
|x|, & |x|\ge2.
\endcases
$$
Note that both expressions agree at $x=\pm2$, so that $\Om(x)$ is
continuous. Moreover, the first derivative $\Om'(x)$ is continuous on
the whole $\R$, while $\Om''(x)$ is not. This is clear from the
explicit expressions
$$
\Om'(x)=\frac2\pi\,\arcsin\frac x2\,, \quad
\Om''(x)=\frac2\pi\,\frac1{\sqrt{4-x^2}}\,, \qquad |x|<2.
$$
Next, we have $|\Om'(x)|<1$ for $|x|<2$, which implies that $\Om(x)$
belongs to $\Dz$. Applying \tht{2.2} for $\om=\Om$ we get
$$
\tp_k[\Om]=-k\int_\R x^{k-1}\left(\frac{\Om(x)-|x|}2\right)'dx
=-\tfrac k2\int_{-2}^2 x^{k-1}(\tfrac2\pi\arcsin\tfrac x2 -\sgn x)dx.
$$

\proclaim{Proposition 5.3} We have
$$
\tp_k[\Om]=\cases \dfrac{(2m)!}{m!m!}\,, & k=2m, \quad m=1,2,\dots, \\
0, & k=1,3,5, \dots\,.
\endcases
$$
\endproclaim

\demo{Proof} Since $\Om(x)$ is even, $\tp_k[\Om]=0$ for odd $k$.  For
even $k$ we get from \tht{2.2}
$$
\tp_{2m}[\Om]=\int_0^2
(1-\tfrac2\pi\,\arcsin\tfrac x2)\,d(x^{2m}).
$$
Setting $x=2\sin\theta$ and integrating by parts we get the result. \qed
\enddemo

We proceed to the ``law of large numbers'' for the Plancherel
measures $M_n$. Actually, it is implied by the ``central limit
theorem'' which will be established in \S7. However, we prefer to
give here an independent short proof.

Recall that in Definition 2.3 we have attached to any Young diagram
$\la\in\Y_n$ a function $\la(\,\cdot\,)\in\Dz$. We define now a
scaled version of it:
$$
\bar\la(x)=n^{-1/2}\la(n^{1/2}x), \qquad x\in\R, \quad \la\in\Y_n\,.
$$
This is a special case of Definition 2.10.  The correspondence
$\la\mapsto\bar\la(\,\cdot \,)$ provides an embedding
$\Y_n\hookrightarrow\Dz$.

\proclaim{Theorem 5.4 (Law of large numbers, 1st form)} Let $\la$ range over
$\Y_n$, and let us view $\bar\la(\,\cdot \,)$ as a random function defined on
the probability space $(\Y_n, M_n)$, where $M_n$ is the Plancherel measure. Let
$\Om$ be as above. Then we have
$$
\lim_{n\to\infty}\int(\bar\la(x)-\Om(x))x^kdx=0 \quad
\text{\rm in probability, for any $k=0,1,\dots$}.\tag5.2
$$
\endproclaim

\demo{Proof} Let $\overline M_n$ be the pushforward of $M_n$ under the
embedding $\Y_n\hookrightarrow\Dz$. Then $\overline M_n$ is a probability
measure on the space $\Dz$.
Given a ``test'' function $f$ on $\Dz$, let $\langle f, \bar
M_n\rangle$ denote the result of pairing between $f$ and $\overline M_n$:
$$
\langle f,\overline M_n\rangle=
\sum_{\la\in\Y_n}f(\bar\la(\,\cdot\,))M_n(\la).
$$
Recall that the elements of $\Alg$ can be interpreted as functions on
$\Dz$, see Definition 2.10 b). Let us take them as ``test''
functions.

We claim that
$$
\lim_{n\to\infty}\langle f, \overline M_n\rangle =f[\Om], \qquad f\in\Alg.
\tag5.3
$$

In other words, the measures $\overline M_n$ on the space $\Dz$ converge
to the Dirac measure at $\Om\in\Dz$ in the weak
topology defined by the function algebra $\Alg$.

Let us prove \tht{5.3}. Without loss of generality we may assume that
$f$ is a homogeneous element with respect to the weight grading in $\Alg$.
Then, by virtue of Proposition 2.11,
$$
\langle f,\overline M_n\rangle=n^{-\wt(f)/2}\,\m f\,. \tag5.4
$$
By Proposition 5.2, $\m f$ is a polynomial in $n$ of degree less
or equal to $\tfrac12\deg_1(f)$. Note that
$\deg_1(\,\cdot\,)\le\deg_\N(\,\cdot\,)$ and $\deg_\N(\,\cdot\,)$
coincides with $\wt(\,\cdot\,)$, see Proposition 4.10. Thus, the
degree of $\m f$ does not exceed $\tfrac12\wt(f)$, which implies that
\tht{5.4} has a limit as $n\to\infty$.

Expand $f$ in the basis $\{\pp_\rho\}$:
$$
f=\sum_\rho f_\rho\pp_\rho\,.
$$
By virtue of Proposition 5.1,
$$
\lim_{n\to\infty} n^{-\wt(f)/2}\,\m f
=f_{(1^{\wt(f)/2})}
$$
with the understanding that the symbol $f_{(1^{k/2})}$ means 0 whenever
$k$ is odd.

Next, Proposition 4.9 (together with Proposition 4.10) implies the
multiplicativity property
$$
(fg)_{(1^{\wt(fg)/2})}\,=f_{(1^{\wt(f)/2})}g_{(1^{\wt(g)/2})}\,,
$$
where $f,g$ are arbitrary weight homogeneous elements. Consequently,
it suffices to examine the case $f=\tp_k$\,, i.e., to show that
$$
(\tp_k)_{(1^{k/2})}\,=\tp_k[\Om].
$$
The right--hand side was found in Proposition 5.3, while the
left--hand side can be evaluated using Proposition 3.7. The result is
the same, which concludes the proof of \tht{5.3}.

Now let us show that \tht{5.3} implies \tht{5.2}. Indeed, \tht{5.2} is
equivalent to
$$
\lim_{n\to\infty}\int\frac{\bar\la(x)-|x|}2\,x^kdx=
\int\frac{\Om(x)-|x|}2x^kdx
\quad\text{in probability, for $k=0,1,\dots$\,.}
$$
By Proposition 2.2, this is equivalent to
$$
\lim_{n\to\infty}\tp_k[\bar\la(\,\cdot\,)]=\tp_k[\Om]
\quad\text{\rm in probability, for $k=2,3,\dots$\,.}
$$
Applying Chebyshev's inequality we see that to prove this, it
suffices to check that the first and the second moments of the random
variable $\tp_k[\bar\la(\,\cdot\,)]$ converge, as $n\to\infty$, to
$\tp_k[\Om]$ and $\tp_k^{\,2}[\Om]$, respectively. But this is a
particular case of \tht{5.3} corresponding to $f=\tp_k$ and
$f=\tp_k^{\,2}$, respectively. \qed
\enddemo

\proclaim{Theorem 5.5 (Law of large numbers, 2nd form)} Let $\la$ range over
$\Y_n$, and let us view $\bar\la(\,\cdot \,)$ as a random function defined on
the probability space $(\Y_n, M_n)$, where $M_n$ is the Plancherel measure. Let
$\Om$ be as above. Then we have
$$
\lim_{n\to\infty}\sup_{x\in\R}|\bar\la(x)-\Om(x)|=0 \quad
\text{\rm in probability.}\tag5.5
$$
\endproclaim

We need two lemmas.

\proclaim{Lemma 5.6} There exists an interval $I\subset\R$ such that the
probability that $\bar\la(x)-|x|$ is supported by $I$ tends to 1 as
$n\to\infty$.
\endproclaim

\demo{Proof} This follows from a finer result due to Hammersley
\cite{Ha}. He has proved that there exists a constant $c$ such that for
any $\varepsilon>0$
$$
\lim_{n\to\infty}M_n\{\la\in\Y_n\mid |\la_1-c\sqrt n|<\varepsilon, \quad
|\la'_1-c\sqrt n|<\varepsilon\}=1.
$$
Actually, the constant $c$ equals 2 (this was first proved by
Vershik--Kerov \cite{VeK1}, \cite{VeK3}), and at present much more is
known about the asymptotics of $\la_1$, see, e.g., the expository
paper \cite{AD}. But, for our purpose, the old Hammersley's result is
enough. \qed
\enddemo

\proclaim{Lemma 5.7} Fix an interval $I=[a,b]\subset\R$, and let $\Sigma$
denote the set of all real--valued functions $\si(x)$ on $\R$, supported by $I$
and satisfying the Lipschitz condition $|\si(x_1)-\si(x_2)|\le|x_1-x_2|$.

On the set $\Sigma$, the weak topology defined by the functionals
$$
\si\mapsto\int\si(x)x^kdx, \qquad k=0,1,\dots,
$$
coincides with the uniform topology defined by the supremum norm
$\Vert\si\Vert=\sup|\si(x)|$.
\endproclaim

This fact was pointed out in \cite{Ke2, \S2.5}.

\demo{Proof} Clearly, the uniform topology is stronger than the weak
topology. Let us check the inverse claim. Given $x\in I$ and
$\varepsilon>0$, let
$$
V(x,\varepsilon)=\{\si\in\Sigma\mid |\si(x)|\le\varepsilon\}.
$$
Pick points $a=\varepsilon_1<\dots<\varepsilon_n=b$ dividing $I$ into
subintervals of length $\le2\varepsilon$. Then, by the Lipschitz
condition, the ball $\Vert\si\Vert\le2\varepsilon$ contains the
intersections of $V(x_i,\varepsilon)$'s. Hence the required claim
reduces to the following one:

Fix $x\in I$ and $\varepsilon>0$. Then $V(x,\varepsilon)$ contains a
neighborhood of 0 in the weak topology.

Let us remark that functions $\si\in\Sigma$ are uniformly bounded,
$\Vert\si\Vert\le(b-a)/2$. Hence the weak topology on $\Sigma$ will
not change if we take, as functionals, integrals with arbitrary
continuous functions $F$. Now let us take a continuous function $F(x)\ge0$,
concentrated in the $\varepsilon/2$--neighborhood of $x$ and such
that $\int F(y)dy=1$. We claim that
$$
\left|\int\si(x)F(x)dx\right|\le\varepsilon/2\, \Rightarrow
\, \si\in V(x,\varepsilon).
$$

Indeed, assume that $\si\notin V(x,\varepsilon)$, i.e.,
$|\si(x)|>\varepsilon$. Without loss of generality we may assume that
$\si(x)>\varepsilon$. Then, for any $y$ such that
$|x-y|\le\varepsilon/2$, we have $\si(y)>\varepsilon/2$, hence
$\int\si(y)F(y)>\varepsilon/2$, which proves our claim. \qed
\enddemo

\demo{Proof of Theorem 5.5} This immediately follows from Theorem 5.4 and
Lemmas 5.6, 5.7. \qed
\enddemo

\head \S6. The central limit theorem for characters \endhead

For any $f\in\Alg$, we denote by $f^{(n)}$ the random variable
defined on the probability space $(\Y_n\,, M_n)$ and obtained by
restricting $f$ to $\Y_n$.

By the symbol $\dc$ we will denote convergence of random variables in
distribution, see, e.g. \cite{Sh}.

The aim of this section is to prove the following result.

\proclaim{Theorem 6.1 (Central limit theorem for characters)}

Choose a sequence $\{\xi_k\}_{k=2,3,\dots}$ of independent standard
Gaussian random variables. As $n\to\infty$, we have
$$
\left\{\frac{{\pp_k}^{(n)}}{n^{k/2}}\right\}_{k=2,3,\dots}\, \dc\,
\{\sqrt k\xi_k\}_{k=2,3,\dots} \tag6.1
$$
\endproclaim

In more detail, for any fixed $N=2,3,\dots$, the joint distribution
of $N-1$ random variables
$$
\frac{{\pp_k}^{(n)}}{\sqrt{k}\,n^{k/2}}\,, \quad 2\le k\le N,
$$
weakly tends, as $n\to\infty$, to the standard Gaussian measure on
$\R^{N-1}$. Note that we could take equally well in \tht{6.1}
the random variables
$$
\frac{n^{k/2}}{\sqrt k}\,\frac{\chi^\la_{(k,1^{n-k})}}{\dim\la}\,,
\qquad 2\le k\le N,
$$
where $\la\in\Y_n$ is the random Plancherel diagram.

The proof of Theorem 6.1 will be given after some preparation work,
based on Propositions 4.11, Proposition 4.12, and Corollary 4.13.

It will be convenient to extend the algebra $\Alg$: we adjoin to it
the square root of the element $\pp_1=p_1$ and then localize over the
multiplicative family generated by $\sqrt{\pp_1}$. Let $\Ae$ denote
the resulting algebra. As a basis in $\Ae$ one can take the elements
of the form
$$
\pp_\rho\cdot(\pp_1)^{m/2}\,, \qquad m_1(\rho)=0, \quad m\in\Z. \tag6.2
$$
We equip $\Ae$ with a filtration by assigning to
$\pp_\rho\cdot(\pp_1)^{m/2}$ the degree
$\deg_1(\,\cdot\,)=|\rho|_1+m$. That is, the
$N$th term of the filtration is spanned by all basis elements \tht{6.2}
with $|\rho|_1+m\le N$. Here $N$ ranges over $\Z$. On the subalgebra
$\Alg\subset\Ae$, this filtration agrees with that induced by the
Kerov degree. Indeed, this claim follows from Proposition 4.11.

Since ${\pp_1}^{(n)}\equiv n$, the symbol $f^{(n)}$ makes sense for
any $f\in\Ae$. Specifically, if $f=g(\pp_1)^{m/2}$ with $g\in\Alg$
and $m\in\Z$ then $f^{(n)}=g^{(n)}\cdot n^{m/2}$. Note also that
Proposition 5.2 admits the following extension:

\proclaim{Proposition 6.2} For any $f\in\Ae$, $\m f$ is a Laurent
polynomial in $n^{1/2}$ whose degree with respect to $n$ is bounded
from above by $\tfrac12\deg_1(f)$.
\endproclaim

\demo{Proof} Evident from Proposition 5.2 and the fact that
${\pp_1}^{(n)}\equiv n$. \qed
\enddemo

Let $H_m(x)$, where $m=0,1,2,\dots$, be the Hermite polynomials in
the normalization of \cite{Sz}, \cite{Er}. We shall need
slightly modified polynomials, which we denote by $\Cal H_m(x)$:
$$
\Cal H_m(x)=2^{-m/2}\, H_m(\sqrt 2x)=
m!\,\sum_{j=0}^{[m/2]}\frac{(-1/2)^jx^{m-2j}}{j!(m-2j)!}\,.
$$
These are monic polynomials, which form the orthogonal system with
respect to standard Gaussian measure $(2\pi)^{-1/2}\exp(-x^2/2)dx$.
They are characterized by the recurrence relation
$$
x\Cal H_m=\Cal H_{m+1}+m\,\Cal H_{m-1} \tag6.3
$$
together with the initial data $\Cal H_0=1$, $\Cal H_1=x$.

For an arbitrary partition $\rho$, we define the element
$\eta_\rho\in\Ae$, which is a normalization of $\pp_\rho$:
$$
\eta_\rho=\frac{\pp_\rho}
{(\pp_1)^{m_1(\rho)}\prod_{k\ge2}(k(\pp_1)^k)^{m_k(\rho)/2}}
=\frac{\pp_\rho}{(\pp_1)^{|\rho|_1/2}\prod_{k\ge2}k^{m_k(\rho)/2}}\,.
\tag6.4
$$
Note that $\deg_1(\eta_\rho)=0$.

We abbreviate $\eta_k=\eta_{(k)}$\,. Note that
$$
\eta_k=\frac{\pp_k}{\sqrt k\,(\pp_1)^{k/2}}\,, \quad k=2,3,\dots\,.
\tag6.5
$$

\proclaim{Proposition 6.3} For any partition $\rho$, we have
$$
\eta_\rho=\prod_{k\ge2}\Cal H_{m_k(\rho)}(\eta_k)+\dots \tag6.6
$$
where dots denote a remainder term with $\deg_1(\,\cdot\,)<0$.
\endproclaim

In particular, $\eta_\rho$ does not depend, up to terms of negative
degree, from the value of $m_1(\rho)$.

\demo{Proof} Examine first the particular case $\rho=(k^m)$, where
$k=2,3,\dots$ and $m=1,2,\dots$\,. Then our claim means that
$$
\eta_{(k^m)}=\Cal H_m(\eta_k)+\dots\,.
$$
Taking $\si=(k^m)$ in Proposition 4.12 we get
$$
\pp_{(k^m)}\cdot \pp_k=
\pp_{(k^{m+1})}+k\,m\,\pp_{(k^{m-1},1^k)}+\dots
=\pp_{(k^{m+1})}+k\,m\,\pp_{(k^{m-1})}\cdot(\pp_1)^k+\dots
$$
for any $m\ge1$, where dots mean lower degree terms. This is equivalent to
$$
\eta_{(k^m)}\cdot \eta_k=\eta_{(k^{m+1})}+m\,\eta_{(k^{m-1})}+\dots,
\qquad m\ge1, \tag6.7
$$
where dots mean terms of negative degree. Within these terms, \tht{6.7}
coincides with the recurrence relation \tht{6.3}, which proves \tht{6.6}.

The case of an arbitrary $\rho$ is reduced to the particular case
$\rho=(k^m)$ using Corollary 4.13. \qed
\enddemo

The next claim is a well--known general result. It justifies the
moment method, which is a convenient tool for checking convergence in
distribution.

\proclaim{Proposition 6.4} Let $a^{(n)}$ be a sequence of real random
variables. Assume that $a^{(n)}$ have finite moments of any order,
and the moments converge, as $n\to\infty$, to the respective moments
of a random variable $a$. Finally, assume that $a$ is uniquely
determined by its moments, which holds, e.g., if the characteristic
function of $a$ is analytic.

Then $a^{(n)}\dc a$. Moreover, this claim also holds when the
variables in question take vector values, i.e., when each $a^{(n)}$,
as well as $a$, is a system of random variables.
\endproclaim

\demo{Sketch of proof} Let $P^{(n)}$ denote the distribution of the
random variable $a^{(n)}$ and $P$ be the distribution of $a$. We have
to prove that $P^{(n)}$ weakly converges to $P$ as $n\to\infty$. The
assumption on the moments implies that $\{P^{(n)}\}$ is a tight
family of probability measures on $\R$. So, it suffices to prove that
any partial weak limit $P'$ of the sequence $\{P^{(n)}\}$ coincides
with $P$. Using again the condition on the moments one can show that
the moments of $P'$ exist and coincide with the limits of the
respective moments of $\{P^{(n)}\}$. Hence, these are exactly the
moments of $P$. By the uniqueness assumption, $P'=P$. \qed
\enddemo

For another proof, see Feller \cite{Fe, ch. VIII, \S6, Example b}.

\demo{Proof of Theorem 6.1} We must prove that
$$
\{\eta_k^{(n)}\}_{k\ge2}\dc \{\xi_k\}_{k\ge2}\,, \qquad n\to\infty.
$$
By Proposition 6.4, it suffices to check that
$$
\m{\prod_{k\ge2}\eta_k^{m_k}}\to
\langle\prod_{k\ge2}\xi_k^{m_k}\rangle_{\operatorname{Gauss}}\,\,,
\qquad n\to\infty, \tag6.8
$$
for any finite collection $\{m_k\}_{k\ge2}$ of nonnegative integers,
where the brackets \linebreak
$\langle\,\cdot\,\rangle_{\operatorname{Gauss}}$
mean expectation with respect to the standard Gaussian measure. The
uniqueness hypothesis of Proposition 6.4 is clearly satisfied.

The limit relations \tht{6.8} are equivalent to the following ones:
$$
\m{\prod_{k\ge2}\Cal H_{m_k}(\eta_k)} \to
\prod_{k\ge2}\langle\Cal H_{m_k}(\xi_k)\rangle_{\operatorname{Gauss}}
$$
for any finite collection $\{m_k\in\Z_+\}$.

If all the numbers $m_k$ are equal to 0 then the expressions in both
sides equal 1, and there is nothing to prove.  So, let us assume that
some of the $m_k$'s are nontrivial. Then the right--hand side
vanishes, because, for a standard Gaussian $\xi$,
$$
\langle\Cal H_m(\xi)\rangle_{\operatorname{Gauss}}
=\tfrac1{\sqrt{2\pi}}\,\int_\R \Cal H_m(x)e^{-x^2/2}dx=0, \qquad
m=1,2,\dots,
$$
by the orthogonal property of the polynomials $\Cal H_m$.

Let us examine the left--hand side. Set
$\rho=(\prod_{k\ge2}k^{m_k})$. By Proposition 6.2,
$$
\prod_{k\ge2}\Cal H_{m_k}(\eta_k)=\eta_\rho\,+\,
\text{a ``remainder term'' of strictly negative degree}.
$$
By our assumption, $\rho$ is nonempty. Moreover, $m_1(\rho)=0$, so
that $\rho\ne(1^r)$. By Proposition 5.1, $\m{\eta_\rho}\equiv0$.
Finally, by Proposition 6.2,
$$
\m{\text{\rm the ``remainder term''}}=O(n^{-1/2}). \tag6.9
$$
This concludes the proof. \qed
\enddemo

Theorem 6.1 can be generalized as follows:

\proclaim{Theorem 6.5} Let $\rho$ range over the set of all
partitions. We have
$$
\{\eta_\rho^{(n)}\}\dc \{\prod_{k\ge2}\Cal H_{m_k(\rho)}(\xi_k)\}\,,
\quad n\to\infty,
$$
where, as before, $\xi_2,\xi_3,\dots$ are independent standard
Gaussians.
\endproclaim

\demo{Proof} The above argument shows that any mixed moment of the
random variables from the left--hand side converges, as $n\to\infty$,
to the respective moment of the random variables from the
right--hand side. However, we cannot use the moment method, because
a polynomial in Gaussian variables does not necessarily satisfy the
uniqueness assumption mentioned in Proposition 6.4. For this reason
we argue in a different way.

Assume that
$$
\{a^{(n)}_1,a^{(n)}_2,\dots\}\dc \{a_1,a_2,\dots\}\,, \quad
n\to\infty,
$$
where $a^{(n)}=\{a^{(n)}_1,a^{(n)}_2,\dots\}$ are families of random
variables depending on $n$ and $a=\{a_1,a_2,\dots\}$ is one more
family of random variables. Next, assume that $f_1(x)=f_1(x_1,x_2,\dots)$,
$f_2(x)=f_2(x_1,x_2,\dots)$, \dots are continuous functions in real
variables $x=(x_1,x_2,\dots)$, where each function actually depends on
finitely many variables only. Then
$$
\{f_1(a^{(n)}),f_2(a^{(n)}),\dots\}\dc
\{f_1(a),f_2(a),\dots\}\,, \quad n\to\infty.
$$

Using this general fact we conclude from Theorem 6.1 that any
polynomial in $\eta^{(n)}_2, \eta^{(n)}_3,\dots$ converges in
distribution to the same polynomial in $\xi_2,\xi_3,\dots$. Moreover,
this also holds for any finite system of polynomials. By virtue
of \tht{6.9}, each $\eta_\rho$ is a polynomial in $\eta^{(n)}_2,
\eta^{(n)}_3,\dots$, within a ``remainder term''. So we only need to
check that the ``remainder term'' does not affect the convergence in
distribution.

Remark that the ``remainder term'' is of the form
$r^{(n)}$, where $r$ is an element of $\Ae$ of strictly negative
degree. It follows that any moment of $r^{(n)}$ tends to 0 as
$n\to\infty$, which implies that $r^{(n)}\dc0$. This shows that the
``remainder term'' is negligible.  \qed
\enddemo

For a different proof of Theorems 6.1 and 6.5, see \cite{Ho}.

\head \S7. The central limit theorem for Young diagrams \endhead

Given $\la\in\Y_n$\,, we set
$$
\De_\la(x)=\tfrac{\sqrt n}2(\bar\la(x)-\Om(x)), \qquad x\in\R. \tag7.1
$$
This is a continuous function on $\R$ with compact support. Dropping
$\la$, which we consider as the random element from the probability
space $(\Y_n,M_n)$, we interpret \tht{7.1} as a random function
$\De^{(n)}(x)$.

For any polynomial $v\in\R[x]$, the integral
$$
v^{(n)}=\int_\R v(x)\De^{(n)}(x)dx \tag7.2
$$
makes sense (because $\De^{(n)}$ is compactly supported) and is a
random variable. We aim to show that the random variables \tht{7.2},
where $v$ ranges over $\R[x]$, are asymptotically Gaussian.

The result will be stated in terms of the Chebyshev
polynomials of the second kind. Instead of the conventional
polynomials $U_k(x)$ (see \cite{Sz}, \cite{Er}) we prefer to
deal with slightly modified polynomials
$$
u_k(x)=U_k(x/2)=\sum_{j=0}^{[k/2]}(-1)^j\binom{k-j}{j}x^{k-2j}\,,
\qquad k=0,1,2,\dots\,. \tag7.3
$$
Note that
$$
u_k(2\cos\theta)=\frac{\sin((k+1)\theta)}{\sin\theta} \tag7.4
$$
and
$$
\int_{-2}^2u_k(x)u_l(x)\frac{\sqrt{4-x^2}}{2\pi}=\delta_{k,l}\,,
\qquad k,l=0,1,2,\dots\,. \tag7.5
$$

\proclaim{Theorem 7.1 (Central limit theorem for Young diagrams)}

According to \tht{7.2}, let
$$
u^{(n)}_k=\int_\R u_k(x)\Delta^{(n)}(dx), \qquad k=1,2,\dots,
$$
and let, as before, $\xi_2,\xi_3,\dots$ stand for a system of independent
standard Gaussians.

We have
$$
\left\{u^{(n)}_k\right\}_{k\ge1}\,\dc \,
\left\{\frac{\xi_{k+1}}{\sqrt{k+1}}\right\}_{k\ge1}\,, \qquad n\to\infty.
$$
\endproclaim

Recall that ``$\dc$'' means convergence in distribution.

Note that $u^{(n)}_0\equiv0$, which explains why we start with $k=1$, not
$k=0$. The theorem is proved at the end of the section. The scheme of the proof
is as follows. We remark that the moments of $\De^{(n)}$ (i.e., the random
variables $v^{(n)}$, where the $v$'s are monomials) are expressed in terms of
the elements $\tp_k$, appropriately centered and scaled. To evaluate the
asymptotics of the corresponding random variables we employ Theorem 6.1. The
main work reduces to expressing the (centered and scaled) elements $\tp_k$
through the elements $\eta_k$ and vice versa, up to lower degree terms.

As in \S6, it is convenient to deal with the extended algebra
$\Ae\supset\Alg$.  We extend the definition of $\deg_1(\,\cdot\,)$ to
$\Ae$ as explained in \S6.

We introduce the elements $q_1,q_2,\ldots\in\Ae$, which are centered
and scaled versions of the elements $\tp_2,\tp_3,\dots$:
$$
q_k=\cases \dfrac{\tp_{k+1}-\frac{(2m)!}{m!m!}(\pp_1)^m}
{(k+1)(\pp_1)^{k/2}}\,,
&\text{if $k$ is odd, $k=2m-1$, where $m=1,2,\dots$,}\\
\dfrac{\tp_{k+1}}{(k+1)(\pp_1)^{k/2}}\,,
&\text{if $k$ is even, $k=2m$, where $m=1,2,\dots$\,.}
\endcases \tag7.6
$$

Since $\tfrac12\tp_2=p_1$ and $p_1=\pp_1$, we have $q_1=0$.

\proclaim{Proposition 7.2} For any $\la\in\Y$,
$$
\int_\R x^k\De_\la(x)dx=\frac{q_{k+1}(\la)}{k+1}\,, \qquad
k=0,1,\dots\,. \tag7.7
$$
\endproclaim

\demo{Proof} Set $n=|\la|$. By the definition of $\De_\la(x)$, see
\tht{7.1},
$$
\int_\R x^k\De_\la(x)dx
=\sqrt n\int_\R x^k\,\frac{\bar\la(x)-|x|}2 dx\,
-\, \sqrt n\int_\R x^k\,\frac{\Om(x)-|x|}2 dx.
$$
By Proposition 2.2, for any $k=0,1,\dots$,
$$
\gather
\sqrt n\int_\R x^k\,\frac{\bar\la(x)-|x|}2 dx
=n^{-(k+1)/2}\int_\R x^k\,\frac{\la(x)-|x|}2 dx\\
=\tfrac1{(k+1)(k+2)}\, \frac{\tp_{k+2}(\la)}{n^{(k+1)/2}}
=\tfrac1{(k+1)(k+2)}\, \frac{\tp_{k+2}}{(\pp_1)^{(k+1)/2}}\,(\la).
\endgather
$$

By Propositions 2.2 and 5.3, for any $k=0,1,\dots$,
$$
\gather
\sqrt n\int_\R x^k\,\frac{\Om(x)-|x|}2 dx
=\frac{\sqrt n\,\tp_{k+2}[\Om]}{(k+1)(k+2)}\\
=\frac{\sqrt n}{(k+1)(k+2)}\,\cdot\,
\cases \dfrac{(2m)!}{m!m!}\,, &\text{if $k$ is even, $k=2m-2$, where
$m=1,2,\dots$,}\\
0, &\text{if $k$ is odd.}
\endcases
\endgather
$$

Combining this with the definition of $q_1, q_2,\dots$, we get
\tht{7.7}. \qed
\enddemo

In order to apply Theorem 6.1 we need the expression of $\pp_k$ in
terms of $q_2,q_3,\dots$ within lower degree terms. We obtain this in
two steps. First, using a trick, we deduce from Proposition 3.7 a formula
expressing any $q_k$ through $\pp_2, \pp_3, \dots$\,, up to lower
degree terms. See Proposition 7.3. Next, we invert this formula, see
Proposition 7.4. One could derive the result directly from
Proposition 3.3 but this way turns out to be more difficult.

\proclaim{Proposition 7.3} For any $k=2,3,\dots$,
$$
q_k=\sum_{j=0}^{[\frac{k-2}2]}\binom kj\,
\frac{\pp_{k-2j}}{(\pp_1)^{(k-2j)/2}}\, +\,\dots, \tag7.8
$$
where dots mean a remainder term with $\deg_1(\,\cdot\,)<0$.
\endproclaim

Note that the elements occurring in the numerator of the right--hand
side are $\pp_2,\pp_3,\dots$ but not $\pp_1$.

\demo{Proof} The claim of the proposition is equivalent to the
following: for any $k=3,4,\dots$,
$$
\gathered
\tp_k=\sum_{j=0}^{[\frac{k-3}2]}\frac{k^{\fd j+1}}{j!}\,
\pp_{k-1-2j}\,(\pp_1)^j
+\cases \dfrac{(2m)!}{m!m!}\,(\pp_1)^m, &\text{\rm if $k$ is even,
$k=2m$,}\\
0, &\text{\rm if $k$ is odd.}
\endcases\\
+\text{\rm terms with $\deg_1(\,\cdot\,)<k-1$.}
\endgathered \tag7.9
$$

We shall deduce this from Proposition 3.7, which expresses $\tp_k$ as
a polynomial in $\pp_1,\pp_2,\dots$, up to terms of lower weight. A
nontrivial point is how to switch from the weight filtration to the
Kerov filtration.

Write the exact expansion of $\tp_k$ through $\pp_1,\pp_2,\dots$,
$$
\tp_k=\sum_\nu a_\nu\,(\pp_1)^{\nu_1}(\pp_2)^{\nu_2}\dots,
$$
where $a_\nu$ are certain coefficients. Let us set
$$
\Vert\nu\Vert=2\nu_1+3\nu_2+4\nu_3+\dots, \quad
\Vert\nu\Vert'=2\nu_1+2\nu_2+3\nu_3+\dots,
$$
so that
$$
\Vert\nu\Vert'=\Vert\nu\Vert-(\nu_2+\nu_3+\dots). \tag7.10
$$

We have
$$
\wt((\pp_1)^{\nu_1}(\pp_2)^{\nu_2}\dots)=\Vert\nu\Vert, \quad
\deg_1((\pp_1)^{\nu_1}(\pp_2)^{\nu_2}\dots)\le\Vert\nu\Vert'.
$$
By Proposition 3.7, we know all coefficients $a_\nu$ with the maximal
value of $\Vert\nu\Vert$ (it equals $k$), while we need all
coefficients with $\Vert\nu\Vert'\ge k-1$.
By \tht{7.10}, $\Vert\nu\Vert'$ does not exceed $k$, and there are 3
possible cases:

$\bullet$ $\Vert\nu\Vert'=\Vert\nu\Vert=k$. Then
$\nu_2=\nu_3=\dots=0$, $k=2\nu_1$, so that $k$ is even; write it as
$2m$. The corresponding monomial is $(\pp_1)^{k/2}$.

$\bullet$ $\Vert\nu\Vert'=k-1$, $\Vert\nu\Vert=k$. Then
$\nu_2+\nu_3+\dots=1$, i.e., exactly one of the numbers
$\nu_2$,$\nu_3$,\dots equals 1. The corresponding monomial is of the
form
$$
\pp_{k-1-2j}(\pp_1)^j\,, \qquad j=\nu_1=0,1,\dots,[\tfrac{k-3}2].
\tag7.11
$$

$\bullet$ $\Vert\nu\Vert'=\Vert\nu\Vert=k-1$. Then
$\nu_2=\nu_3=\dots=0$, $k$ is odd, and the corresponding monomial is
$(\pp_1)^{(k-1)/2}$.

In the first and second cases, $\Vert\nu\Vert$ takes the maximal
value $k$, and then the coefficients $a_\nu$ are known from
Proposition 3.7. In the third case, $\Vert\nu\Vert$ is no longer
maximal, so that Proposition 3.7 does not tell us what is the
coefficient. However, an additional argument will imply that it is
actually 0.

Indeed, in the third case $k$ must be odd, which implies that
$\tp_k$ is antisymmetric with respect to ``$\inv$'', see \tht{2.17}. On
the other hand, by virtue of \tht{4.10}, any monomial in
$\pp_1,\pp_2,\dots$ is either symmetric or antisymmetric. Therefore,
in the expansion of $\tp_k$ only antisymmetric monomials can occur.
Since $\pp_1$ is symmetric, the monomial $(\pp_1)^{(k-1)/2}$ is also
symmetric, so that it does not appear.

Thus, we have proved that the top degree component of $\tp_k$ is
obtained from the top weight terms as given in Proposition 3.7; we
simply keep all terms proportional either to $(\pp_1)^{k/2}$ or to a
monomial of the form \tht{7.11}, and remove all the remaining terms. This
procedure leads to \tht{7.9}. \qed
\enddemo

In the next proposition we invert \tht{7.8}.

\proclaim{Proposition 7.4} For any $k\ge2$,
$$
\frac{\pp_k}{(\pp_1)^{k/2}}=
\sum_{j=0}^{[\tfrac{k-2}2]}(-1)^j\,\frac{k}{k-j}\,
\binom{k-j}j\, q_{k-2j}\,+\dots, \tag7.12
$$
where dots mean a remainder term with $\deg_1(\,\cdot\,)<0$.
\endproclaim

\demo{Proof} We employ the following combinatorial inversion formula,
see \cite{Ri, \S2.4, (10)}:

Let $a_0,a_1,a_2,\dots,b_0,b_1,b_2,\dots$ be formal variables. Then
$$
\multline
\left\{a_k=\sum_{j=0}^{[k/2]}
\binom kj\, b_{k-2j}\right\}_{k=0,1,\dots}\\
\Longleftrightarrow
\left\{b_k=\sum_{j=0}^{[k/2]} (-1)^j\,\frac k{k-j}\,
\binom{k-j}j\, a_{k-2j}\right\}_{k=0,1,\dots}
\endmultline
\tag7.13
$$

Set
$$
a_0=a_1=0, \quad a_k=q_k \quad (k\ge2); \qquad
b_0=b_1=0, \quad b_k=\frac{\pp_k}{(\pp_1)^{k/2}} \quad (k\ge2).
$$
The relations \tht{7.8} coincide with the first system in \tht{7.13}, up to
remainder terms of negative degree. These terms can be neglected,
because they affect only similar remainder terms in the inverse
relations. These inverse relations are given then by the second system
in \tht{7.13}. This leads to \tht{7.12}. \qed
\enddemo

\demo{Proof of Theorem 7.1} We rewrite \tht{7.12} as follows. For any
$k\ge2$,
$$
\sum_{j=0}^{[\tfrac{k-2}2]}(-1)^j\,\binom{k-1-j}j\,\frac{q_{k-2j}}{k-2j}
=\frac1k\,\frac{\pp_k}{(\pp_1)^{k/2}}\,+R_k\,,
$$
where $R_k\in\Ae$ is a certain element such that $\deg_1(R_k)<0$.

In the left--hand side, we may extend the summation up to
$[\frac{k-1}2]$, because $q_1=0$. Comparing this with formula \tht{7.3}
for $u_{k-1}(x)$ and formula \tht{7.7} for the moments of $\De_\la$, we
conclude that
$$
u^{(n)}_{k-1}=\tfrac1{\sqrt k}\eta^{(n)}_k\,+R^{(n)}_k\,
\qquad k\ge2.
$$
Or, equivalently,
$$
u^{(n)}_k=\tfrac1{\sqrt{k+1}}\eta^{(n)}_{k+1}\,+R^{(n)}_{k+1}\,
\qquad k\ge1.
$$

As $n\to\infty$, the asymptotics of the (mixed) moments of the random
variables $u^{(n)}_1\,$, $u^{(n)}_2\,, \dots$ is the same as that
for the random variables $\tfrac1{\sqrt{k+1}}\eta^{(n)}_{k+1}$\,,
$k=1,2,\dots$\,. Indeed, the remainder terms of negative degree do
not affect the asymptotics, see Proposition 6.2. As for the moments
of the random variables $\eta^{(n)}_{k+1}$, their asymptotics has
been evaluated in the proof of Theorem 6.1. This concludes the proof.
\qed
\enddemo

\head \S8. The central limit theorem for transition measures of Young
diagrams \endhead

Let $\Cal M$ denote the set of
probability measures on $\R$ with compact support, and let
$\Mz\subset\Cal M$ be the subset of measures with the first moment equal
to 0.

\proclaim{Proposition 8.1}
There exists a bijective correspondence $\om\leftrightarrow\mu$
between $\Cal D$ and $\Cal M$, which is also a bijection
$\Dz\leftrightarrow\Mz$. It is characterized by the relation
$$
\exp\int_\R\frac{\si'(x)dx}{x-z}\,=
\int_\R \frac{\mu(dx)}{1-\frac xz}\,, \tag8.1
$$
where $z\in\C\setminus I$, where $I\subset\R$ stands for a
sufficiently large interval, and, as usual,
$\si(x)=\tfrac12(\om(x)-|x|)$.
\endproclaim

\demo{Proof} See \cite{Ke2}, \cite{Ke4}. \qed
\enddemo

We call $\mu $ the {\it transition measure\/} of the continual diagram
$\om$. In \cite{Ke4}, the correspondence $\si'\mapsto \mu$
defined by \tht{8.1} is defined in a greater generality, so that its
range is the set of all (not necessarily compactly supported)
probability measures on $\R$. (Note that in \cite{Ke4}, the symbol
$\Cal M$ refers to the latter set.)

Formula \tht{8.1} means that the two sequences,
$$
\left\{-k\int_\R x^{k-1}\si'(x)dx\right\}_{k=1,2,\dots}
\quad \text{and} \quad
\left\{\int_\R x^k\mu(dx)\right\}_{k=1,2,\dots}\,,
$$
are related to each other in exactly the same way as the two systems
of generators of the algebra $\La$, $\{\bp_k\}_{k=1,2,\dots}$ and
$\{\bh_k\}_{k=1,2,\dots}$.

{}From now on and up to the end of this section we restrict ourselves
to measures $\mu$ from the subset $\Mz$.

\example{Definition 8.2} Recall that the algebra $\Alg$ can
be realized as the image of $\La$ under the morphism \tht{2.18}, and
let $\th_2, \th_3,\dots$ denote the image in
$\Alg$ of the elements $\bh_2,\bh_3,\dots$\,.

We realize $\Alg$ as an algebra of functions on $\Mz$ by setting
$$
\th_k[\mu]=\int_\R x^k\mu(dx), \qquad
\mu\in\Mz, \quad k=2,3,\dots\,.
$$
Equivalently, for any $f\in\Alg$, we set $f[\mu]=f[\om]$, where
$\om\leftrightarrow\mu$ and $f[\om]$ is understood according to
Definition 2.10 b).
\endexample

\proclaim{Proposition 8.3} Let $\om=\la(\,\cdot\,)$, where $\la$ is a
Young diagram, and let
$\{x_i\}_{i=1}^{m+1}$\,, $\{y_j\}_{j=1}^m$ be the local extrema of
$\la(x)$, see \S2.
Then the transition measure of $\om$ is supported by the finite set
$\{x_1,\dots,x_{m+1}\}$ and is given by the following formula:
$$
\mu=\sum_{i=1}^{m+1}\mu_i\delta_{x_i}, \tag8.2
$$
where $\delta_x$ denotes the Dirac mass at $x$ and the weights
$\mu_1, \dots, \mu_{m+1}$ are the coefficients in
the expansion
$$
\frac{\prod_{j=1}^m(z-y_j)}{\prod_{i=1}^{m+1}(z-x_i)}=
\sum_{i=1}^{m+1}\frac{\mu_i}{z-x_i}\,. \tag8.3
$$
\endproclaim

\demo{Proof} See \cite{Ke2}, \cite{Ke4}. \qed
\enddemo

The measure $\mu$ defined by \tht{8.2}--\tht{8.3} is called the {\it
transition measure\/} of a given Young diagram $\la$. For a
justification of this term and more details, see \cite{Ke2},
\cite{Ke3}, \cite{Ke4}.

\proclaim{Proposition 8.4} The transition measure of $\Om\in\Dz$ is
the ``semi--circle distribution'' $\sc$ supported by $[-2,2]$,
$$
\sc(dx)=\tfrac1{2\pi}\sqrt{4-x^2}. \tag8.4
$$
\endproclaim

\demo{Proof} See \cite{Ke2}, \cite{Ke4}. \qed
\enddemo

\example{Definition 8.5} Fix $n=1,2,\dots$. To any $\la\in\Y_n$ we
assign a probability measure $\hla\in\Mz$ as follows: $\hla$ is the
transition measure of the scaled diagram $\bar\la(\,\cdot\,)\in\Dz$.
Equivalently, $\hla$ is the push--forward of the transition
measure \tht{8.2}--\tht{8.3} under the shrinking $x\mapsto n^{-1/2}x$ of
the real axis.

Viewing $\la$ as the random element of the probability space
$(\Y_n,M_n)$, we interpret $\hla$ as a random probability
measure.
\endexample

The next result is simply a reformulation of Theorem 5.5. Recall
that, by Definition 8.2, we may view $\Alg$ as a function algebra on
$\Mz$.

\proclaim{Theorem 8.6 (Law of large numbers for transition measures)}

As $n\to\infty$, the random measures $\hla$ concentrate near the
Dirac mass at the element $\sc\in\Mz$, the semi--circle distribution
\tht{8.4}.

In more detail, let $\widehat M_n$ stand for the push--forward of the
measure $M_n$ under the correspondence $\la\mapsto\hla$. Then
$$
\lim_{n\to\infty}\langle f, \widehat M_n\rangle = f[\sc],
\qquad \forall f\in\Alg.
$$
\endproclaim

\demo{Proof} Immediately follows from Theorem 5.5, Proposition 8.4
and Definition 8.2. \linebreak \qed
\enddemo

Now, our aim is to describe the fluctuations of the random
measures $\hla$ around the semi--circle distribution $\sc$. We do not
know if this can be achieved by a simple application of Theorem 7.1.
The reason is that the transform $\bar\la(\,\cdot\,)\mapsto\hla$ is
highly nonlinear. It turns out, however, that the proof of Theorem
7.1 can be readily translated to the language of transition measures:
it suffices to deal with $\th_2,\th_3,\dots$ instead of
$\tp_2,\tp_3,\dots$\,.

The role of the polynomials $u_k(x)$ is played now by the polynomials
$t_k(x)$. These are slightly modified Chebyshev polynomials of the
first kind. By definition,
$$
t_k(x)=2T_k(\tfrac x2)
=\sum_{j=0}^{[k/2]}(-1)^j\,\frac k{k-j}\,\binom{k-j}j\,x^{k-2j}
=k\sum_{j=0}^{[k/2]}(-1)^j\,\frac{(k-1-j)!}{j!(k-2j)!}\, x^{k-2j}\,,
\tag8.5
$$
where $k=1,2,\dots$ and the $T_k$'s are the conventional Chebyshev
polynomials, see \cite{Sz}, \cite{Er}.

We also have (cf. \tht{7.4}, \tht{7.5})
$$
t_k(2\cos\theta)=2\cos(k\theta), \qquad k=1,2,\dots, \tag8.6
$$
and
$$
\int_{-2}^2t_k(x)t_l(x)\,\frac1{2\pi\,\sqrt{4-x^2}}\,dx=\delta_{k,l}\,,
\qquad k,l=1,2,\dots\,. \tag8.7
$$

Given $\la\in\Y_n$\,, we set (cf. \tht{7.1})
$$
\widehat\Delta_\la=\sqrt n(\hla-\sc). \tag8.7a
$$
This is a compactly supported measure on $\R$ (in general, not a
positive one). Dropping $\la$, which is viewed as the random element
of $(\Y_n,M_n)$, we interpret \tht{8.7a} as a random measure, which we
denote by $\widehat\Delta^{(n)}$.

Next, we set
$$
t^{(n)}_k=\int_\R t_k(x)\,\widehat\Delta^{(n)}(dx).
$$
This is a random variable, defined on the probability space
$(\Y_n,M_n)$.

\proclaim{Proposition 8.7} We have $t^{(n)}_1=t^{(n)}_2\equiv0$.
\endproclaim

\demo{Proof} Recall that for any measure from $\Mz$, the first moment equals
zero. In particular, this holds for $\hla$ and $\sc$, which implies
$t^{(n)}_1\equiv0$.

Next, the relation $\bh_2=\frac12(\bp_1^2+\bp_2)$ in the algebra
$\La$ turns, under the morphism \tht{2.18}, into the relation
$\th_2=\tp_2/2=p_1$ in the algebra $\Alg$. It follows that
$\th_2[\la(\,\cdot\,)]\equiv n=|\la|$, which implies
$\th_2[\bar\la(\,\cdot\,)]\equiv1$. This in turn means that the
second moment of $\hla$ equals 1. On the other hand, the second moment
of $\sc$ also equals 1. Therefore, the second moment of $\hla-\sc$
equals 0 for any $\la\in\Y_n$, so that $\th_2^{(n)}\equiv0$. \qed
\enddemo

\proclaim{Theorem 8.8 (Central limit theorem for transition
measures, cf. Theorem 7.1)} \footnote{As was already mentioned in
Introduction, this result is due to the authors.}

In the notation introduced above,
$$
\left\{t^{(n)}_k\right\}_{k\ge3} \dc
\left\{\sqrt{k-1}\,\xi_{k-1}\right\}_{k\ge3}\,\,,
$$
where $\xi_2,\xi_3,\dots$ are independent standard Gaussian random
variables.
\endproclaim

Here we start with $k=3$, because $t^{(n)}_1=t^{(n)}_2\equiv0$, see
Proposition 8.5.

\demo{Outline of proof} Since the argument is strictly parallel to
that given above for Theorem 7.1, we will not repeat all the details.

\demo{Step 1: Expressing $\th_k$ through $\pp_1,\pp_2,\dots$, up to
lower weight terms}

This is a counterpart of Proposition 3.7. We start with formula
\tht{3.6} of Proposition 3.7, which we rewrite as follows
$$
\pp_{k-1}=-\tfrac1{k-1}\,[t^k]\{A^{-(k-1)}(t)\}+\dots,
\qquad k=2,3,\dots,
$$
where
$$
A(t)=1+\sum_{j=2}^\infty \th_jt^j
$$
and dots mean lower weight terms. Applying Proposition 3.6 we invert
this formula and get
$$
\th_k=\frac1{k+1}[u^k]\{B^{k+1}(u)\}+\dots,
\qquad k=2,3,\dots,
$$
where
$$
B(u)=1+\sum_{j=2}^\infty \pp_{j-1}u^j
$$
and dots mean a polynomial in $\{\pp_{j-1}\}$ of weight $<k$, where,
by definition, \linebreak $\wt(\pp_{j-1})=j$.  More explicitly,
$$
\th_k=\sum\Sb m_2,m_3,\dots\\ 2m_2+3m_3+\dots=k\endSb
\frac{k^{\fd (\sum m_i -1)}}{\prod m_i!}\,
\prod(\pp_{i-1})^{m_i}\,+\,\text{lower weight terms.} \tag8.8
$$
\enddemo

\demo{Step 2: Switching to the Kerov filtration}

The same argument as that used in the proof of Proposition 7.3
makes it possible to derive from \tht{8.8} the following expression, cf.
\tht{7.9}.
$$
\gather
\th_k=\sum_{j=0}^{[\frac{k-3}2]}
\binom kj\pp_{k-1-2j}(\pp_1)^j \, +\,
\cases \tfrac1{m+1}\binom{2m}m (\pp_1)^m\,, &\text{if
$k=2m$, $m=1,2,\dots$,}\\
0, &\text{if $k$ is odd}
\endcases\\
+\,\text{\rm terms with $\deg_1(\,\cdot\,)<k-1$.}
\endgather
$$

Define the elements $g_k\in\Ae$ by (cf. \tht{7.6})
$$
g_k=\cases \dfrac{\th_k-\tfrac1{m+1}\binom{2m}m\,(\pp_1)^m}
{(\pp_1)^{(k-1)/2}}\,, &\text{if $k=2m$, $m=1,2,\dots$,}\\
\dfrac{\th_k}{(\pp_1)^{(k-1)/2}}\,, &\text{if $k$ is odd}.
\endcases \tag8.9
$$
The above expression for $\th_k$ is equivalent to
$$
g_k=\sum_{j=0}^{[\frac{k-3}2]}\binom kj\,
\frac{\pp_{k-1-2j}}{(\pp_1)^{(k-1-2j)/2}}\,+\dots,
\qquad k=3,4,\dots, \tag8.10
$$
where dots mean terms with $\deg_1(\,\cdot\,)<0$.
\enddemo

\demo{Step 3: Inverse formula expressing
$\pp_{k-1}/(\pp_1)^{(k-1)/2}$ through $g_3,g_4,\dots$}

This is a counterpart of Proposition 7.4. We note that \tht{8.10} is
quite similar to \tht{7.8}. Exactly as in Proposition 7.4, we get
$$
\frac{\pp_{k-1}}{(\pp_1)^{(k-1)/2}}=
\sum_{j=0}^{[k/2]}(-1)^j\,\frac k{k-j}\binom{k-j}j\, g_{k-2j}+\dots,
\tag8.11
$$
where, by convention, $g_0=g_1=g_2=0$ and dots mean terms with
$\deg_1(\,\cdot\,)<0$.
\enddemo

\demo{Step 4: Interpretation in terms of Chebyshev's polynomials}

This final step is similar to the proof of Theorem 7.1 at the end of
\S7. The moments of the semi--circle distribution have the following
form, cf. Proposition 5.3:
$$
\int_{-2}^2 x^k\sc(dx)=\cases \tfrac1{m+1}\binom{2m}m\,, &\text{if
$k=2m$, $m=1,2,\dots$,}\\ 0, &\text{if $k$ is odd.}
\endcases
$$
{}From this, the definition of the elements $g_k$\, (see
\tht{8.9}), and Definition 8.5 we get
$$
g_k(\la)=\sqrt n\int_{-\infty}^\infty x^k(\hla-\sc)(dx)
=\int_{-\infty}^\infty x^k\widehat\Delta_\la(dx), \qquad
k=3,4,\dots\,. \tag8.12
$$

It follows from \tht{8.12} that any polynomial $v\in\R[x]$ may be
identified with an element of $\Ae$ (say, $\widehat v$) via
$$
\widehat v(\la)=\int_\R v(x)\widehat\Delta_\la(dx), \qquad \la\in\Y.
$$
Or, equivalently,
$$
\R[x]\ni v=\sum_i c_ix^i \quad \longleftrightarrow\quad
\widehat v=\sum_i c_ig_i\in\Ae.
$$

In particular, for $v=t_k$ we get from \tht{8.5}
$$
\widehat t_k=\sum_{j=0}^{[k/2]}(-1)^j\,\frac k{k-j}\,
\binom{k-j}j\, g_{k-2j}\,, \qquad k=3,4,\dots\,.
$$
Recall also that $\widehat t_0=\widehat t_1=\widehat t_2=0$.

Comparing this with \tht{8.11} we see that
$$
\widehat t_k=\frac{\pp_{k-1}}{(\pp_1)^{(k-1)/2}}\,+\, R_{k-1}\,,
\qquad R_{k-1}\in\Ae, \quad \deg_1(R_{k-1})<0, \qquad k=3,4, \dots\,.
$$
Or by the definition of the elements $\eta_k$, see \tht{6.5},
$$
\widehat t_k=\sqrt{k-1} \eta_{k-1}+R_{k-1}, \qquad k=3,4,\dots,
$$
so that
$$
t^{(n)}_k=\sqrt{k-1}\eta^{(n)}_{k-1}+R^{(n)}_{k-1}\,,
\qquad k=3,4,\dots\,.
$$
Then the proof is completed as at the end of \S7. \qed
\enddemo

\enddemo

\head \S9. Discussion \endhead

Recall that a {\it generalized Gaussian process\/} is a Gaussian
measure in a space $\Cal F'$ of distributions (=generalized functions);
$\Cal F'$ is supposed to be the dual to a
space $\Cal F$ of test functions. Given a test function
$\varphi\in\Cal F$, the
result of its pairing with the random distribution defined by the
process is a random Gaussian variable. See, e.g.,
Gelfand--Vilenkin \cite{GV}, Simon \cite{Si}.

We shall define Gaussian processes via expansions in some orthogonal systems of
functions with random coefficients (a useful general reference on such random
series is Kahane's book \cite{Ka}). Consider the random series
$$
\De(x)=\sum_{k=1}^\infty \frac{\xi_{k+1}u_k(x)\sqrt{4-x^2}}
{2\pi\sqrt{k+1}}\,,
\qquad -2\le x\le2. \tag9.1
$$
Here, as above, $\xi_2$, $\xi_3$,\dots are independent standard Gaussian
random variables and $u_1,u_2,\dots$ are modified Chebyshev's
polynomials of the second kind (see \tht{7.3}, \tht{7.4}). The series
\tht{9.1} correctly defines a generalized Gaussian process, where as
$\Cal F'$ we take the space $(C^\infty(\R))'$ of compactly supported
distributions on the real line. But the process actually lives on the
subspace of distributions concentrated on $[-2,2]$). For any test function
$\varphi\in\Cal F=C^\infty(\R)$,
$$
\langle\varphi,\De\rangle=\int_{-2}^2\varphi(x)\De(x)dx=
\sum_{k=1}^\infty\frac{\xi_{k+1}}{2\pi\sqrt{k+1}}\,
\int_{-2}^2\varphi(x)u_k(x)\sqrt{4-x^2}dx
$$
is a Gaussian random variable. In particular, setting
$\varphi(x)=u_k(x)$ we get, by the orthogonality relation \tht{7.5},
$$
\langle u_k,\De\rangle=\frac{\xi_{k+1}}{\sqrt{k+1}}\,, \qquad
k=1,2,\dots\,.
$$

Informally, the result of Theorem 7.1 can be stated as follows: for
the random Plancherel diagram $\la\in\Y_n$\,,
$$
\bar\la(x)\sim \Om(x)+\frac2{\sqrt n}\,\De(x), \qquad n\to\infty, \tag9.2
$$
where $\bar\la(x)$ was introduced in Definition 2.3, and $\De(x)$ is
given by \tht{9.1}

Next, consider the random series
$$
\widehat\De(x)=\sum_{k=3}^\infty \frac{\sqrt{k-1}\xi_{k-1}t_k(x)}
{2\pi\sqrt{4-x^2}}\,,
\qquad -2<x<2. \tag9.3
$$
Here $\xi_2,\xi_3,\dots$ are as above and $t_3,t_4,\dots$ are
modified Chebyshev's polynomials of the first kind (see \tht{8.5},
\tht{8.6}). The series \tht{9.3} correctly defines a generalized
Gaussian process on the same space $\Cal F'=(C^\infty(\R))'$ of compactly
supported distributions. For any test function
$\varphi\in C^\infty(\R)$,
$$
\langle\varphi,\widehat\De\rangle=\int_{-2}^2\varphi(x)\widehat\De(x)dx=
\sum_{k=3}^\infty\frac{\sqrt{k-1}\xi_{k-1}}{2\pi}\,
\int_{-2}^2\frac{\varphi(x)u_k(x)}{\sqrt{4-x^2}}dx
$$
is a Gaussian random variable. In particular, setting
$\varphi(x)=t_k(x)$ we get, by the orthogonality relation \tht{8.7},
$$
\langle t_k,\widehat\De\rangle=\sqrt{k-1}\xi_{k-1}\,, \qquad
k=3,4,\dots\,.
$$

Informally, the result of Theorem 8.6 can be stated as follows: for
the random Plancherel diagram $\la\in\Y_n$\,,
$$
\widehat\la(x)\sim \sc(x)+\frac1{\sqrt n}\,\widehat\De(x),
\qquad n\to\infty, \tag9.4
$$
where the transition measure $\widehat\la$ (see Definition 8.5) is
viewed as a generalized function, and $\widehat\De(x)$ is given by
\tht{9.3}.

Let us compare these results with the central limit theorem for
the Gaussian unitary ensemble. Consider the space $\Bbb H_N$ of
$N\times N$ complex Hermitian matrices, and equip it with the
Gaussian measure
$$
\Gauss_N(dX)=\left(\frac N{2\pi}\right)^{N^2/2}\,
\exp\left\{-\frac N2\,\operatorname{tr} (X^2)\right\}\Leb(dX),
\tag9.5
$$
where $X$ ranges over $\Bbb H_N$ and $\Leb$ denotes the Lebesgue
measure on $\Bbb H_N\simeq\R^{N^2}$.

To any matrix $X\in\Bbb H_N$ we assign a certain probability measure
$\mu_X$ on $\R$, which we prefer to view as a generalized function:
$$
\mu_X(x)=\frac1N(\de(x-x_1)+\dots+\de(x-x_N)), \tag9.6
$$
where $x_1\,,\dots, x_N$ are the eigenvalues of $X$, and $\de(x)$ is
the delta function. Dropping $X$, which we view as the random element
of the probability space $(\Bbb H_N,\Gauss_N)$, we regard \tht{9.6}
as the random generalized function $\mu(x)$. Then we have
the following {\it central limit theorem for the Gaussian unitary
ensemble\/} (we state it informally):
$$
\mu(x)\sim\sc(x)+\frac1N\,\widetilde\De(x), \qquad N\to\infty, \tag9.7
$$
where $\widetilde\De(x)$ is the Gaussian process on $[-2,2]$ defined
by the random series
$$
\widetilde\De(x)
=\sum_{k=1}^\infty\frac{\sqrt k\,\xi_k t_k(x)}{2\pi\sqrt{4-x^2}} \tag9.8
$$
with independent standard Gaussians $\xi_1,\xi_2,\dots$\,.

For the rigorous formulation and proof of this result (and its
generalizations), see \cite{Jo2}. Note that similar results hold for
other random matrix ensembles, see \cite{DE}, \cite{DS}, \cite{Jo1}.
As explained in \cite{Jo1}, \cite{Jo2}, this subject has close links
with the famous Szeg\"o theorem on asymptotics of Toeplitz determinants.

Comparing \tht{9.3} and \tht{9.8} we see that the Gaussian processes
$\widehat\De(x)$and $\widetilde\De(x)$ look rather close. Another
observation is that
$$
-\frac12\cdot\frac d{dx}\, \De(x)
=\sum_{k=2}^\infty\frac{\sqrt k\,\xi_k t_k(x)}{2\pi\sqrt{4-x^2}}\,.
\tag9.9
$$
That is, the derivative of the process $\De(x)$ coincides, up to
factor $-1/2$ and the first term, with the process
$\widetilde\De(x)$. This is readily seen from the following formulas.

All the three series $\De(x)$, $\widehat\De(x)$, $\widetilde\De(x)$
look especially simply after change of a variable, $x=2\cos\theta$.
Using \tht{7.4}, \tht{8.6} we get
$$
\gather
\De(2\cos\theta)
=\frac1\pi\,\sum_{k=2}^\infty\frac{\xi_k}{\sqrt k}\,\sin(k\theta)\,, \\
\widehat\De(2\cos\theta)
=\frac1{2\pi}\,\sum_{k=3}^\infty\frac{\sqrt{k-1}\,\xi_{k-1}\cos(k\theta)}
{\sin\theta}\,,\\
\widetilde\De(2\cos\theta)
=\frac1{2\pi}\,\sum_{k=1}^\infty\frac{\sqrt k\,\xi_k\cos(k\theta)}
{\sin\theta}\,.
\endgather
$$

\head \S10. Free cumulants and Biane's theorem \endhead

Let, as above,
$$
H(t)=1+\sum_{j=1}^\infty \bh_jt^j\,\in\La[[t]]
$$
be the generating series for the complete homogeneous symmetric
functions. We introduce elements $\ff_1$\,, $\ff_2$\,, \dots in $\La$
as follows:
$$
\ff_1=\bh_1\,; \qquad
\ff_k=-\,\frac1{k-1}\,[t^k]\left\{H^{-(k-1)}(t)\right\},
\qquad k=2,3,\dots\,.
$$
More explicitly,
$$
\ff_k=\sum\Sb m_1,m_2,\dots\ge0\\
1\cdot m_1+2\cdot m_2+\dots=k\endSb
(-1)^{\sum m_j-1}\,k^{\uparrow(\sum m_j-1)}\,
\prod_{j\ge1}\frac{\bh_j^{m_j}}{m_j!}\,,
$$
where
$$
x^{\uparrow m}=\frac{\Gamma(x+m)}{\Gamma(x)}=x(x+1)\dots(x+m-1).
$$

This definition is inspired by Voiculescu's free probability theory
\cite{Vo}, \cite{VoDN}.
Let $\mu$ be a compactly supported probability measure on $\R$, i.e.,
an element of $\Cal M$, in our notation. When the $\bh_k$'s are
specialized to the moments of $\mu$,
$$
\bh_k\quad \longrightarrow\quad \int_\R x^k\mu(dx), \qquad
k=1,2,\dots,
$$
the elements $\ff_k$ turn into the {\it free cumulants\/} of the
measure $\mu$. The free cumulants are
counterparts of the semi--invariants in the sense of conventional
probability theory. The free cumulants are additive functionals with
respect to additive free convolution of measures (just as the
semi--invariants are additive functionals with respect to the
conventional convolution product). See \cite{Vo}, \cite{VoDN}, \cite{Sp}.

Denote by $\tf_k\in\Alg$ the image of $\ff_k\in\La$ under the
morphism \tht{2.18}. Note that $\tf_1=0$.
Let $\la\in\Y$ be arbitrary and let $\mu\in\Cal M^0$ be the
transition measure of $\la$, see \tht{8.2}--\tht{8.3}. Then
$\tf_k(\la)$ coincides with the $k$th free cumulant of $\mu$.

\proclaim{Proposition 10.1} For any $k=1,2,\dots$, the element
$\tf_{k+1}$ coincides with the top weight homogeneous component of
$\pp_k$.
\endproclaim

\demo{Proof} By the very definition
$$
\tf_{k+1}=-\frac1k\,[t^{k+1}]\left\{(1+\sum_{j\ge2}\th_jt^j)^{-k}\right\}
=-\frac1k\,[t^{k+1}]\left\{\exp(-k\sum_{j\ge2}\frac{\tp_j}j\,t^j)\right\}\,,
$$
which is exactly formula \tht{3.4}. Then the claim follows from
Proposition 3.5. \qed
\enddemo

Biane \cite{Bi1} found out that free cumulants emerge in the asymptotic
theory of characters of the symmetric groups. To state his result we
need a notation.

Given $A>1$, let $\Y(A)$ denote the set of the Young diagrams $\la$
such that $\la_1\le\sqrt n A$, $\la'_1\le\sqrt n A$, where $n=|\la|$.
Equivalently, $\bar\la(x)=|x|$ whenever $|x|\ge A$. Recall that
$\bar\la(x)=n^{-1/2}(n^{1/2}x)$ is the scaled version of
$\la(\,\cdot\,)$.

\proclaim{Theorem 10.2 (\cite{Bi1, Th. 1.3})} Fix an arbitrary $A>1$.
For any partition $\rho$ and any $\la\in\Y(A)$,
$$
\frac{\chi^\la_{\rho\cup1^{n-|\rho|}}}{\dim\la}
=n^{-\frac{|\rho|-\ell(\rho)}2}\,
\prod_{j\ge1}\tf_{j+1}^{\,m_j(\rho)}[\bar\la(\,\cdot\,)]\,
+\,O\left(n^{-\frac{|\rho|-\ell(\rho)}2-1}\right),  \tag10.1
$$
where $n=|\la|$ is assumed to be $\ge|\rho|$. Here the estimate of
the remainder term depends only on $A$ and $\rho$, and is uniform on
$\la$ provided that $\la$ ranges over $\Y(A)$.
\endproclaim

\demo{Comments} 1) All terms in \tht{10.1} do not depend on $m_1(\rho)$.
Indeed, this is evident for $\rho\cup1^{n-|\rho|}$ and
$|\rho|-\ell(\rho)$. On the other hand $\tf_2(\bar\la(\,\cdot\,))=1$,
so that the factor corresponding to $j=1$ equals 1.

2) As pointed out by Biane, formula \tht{10.1} implies that if
$\{\la\}$ is a sequence of diagrams in $\Y(A)$ such that
$n=|\la|\to\infty$ and $\bar\la(\,\cdot\,)$ uniformly converges to a
continual diagram $\om\in\Dz$ then
$$
\frac{\chi^\la_{\rho\cup1^{n-|\rho|}}}{\dim\la}\,
\sim \, C\,n^{-\frac{|\rho|-\ell(\rho)}2}\,,
\qquad C=\prod_{j\ge1}\tf_{j+1}^{\,m_j(\rho)}[\om].
\tag10.2
$$
Note that in some cases the constant $C$ can vanish, which implies a faster
decay of the character values: this happens, for instance, when $\om=\Om$ and
$\rho$ is nontrivial (i.e. distinct from $(1^r)$), because $\tf_k[\Om]=0$ for
any $k\ge3$.

3) In this result, the assumption that $\la$ ranges over a set of the form
$\Y(A)$ plays a key role. When this assumption is dropped, quite a different
estimate for the left--hand side of \tht{10.1} holds, see Roichman's paper
\cite{Ro}.

4) Biane \cite{Bi2} obtained further results in this direction.
\enddemo

We shall give an alternative proof of this Biane's theorem. Our
argument seems to be rather simple and transparent.

\demo{Proof of Theorem 10.2} Multiply both sides of \tht{10.1} by
$$
n^{\fd|\rho|}=n^{|\rho|}\, (1+O(n^{-1})).
$$
Then, by Definition 4.1, \tht{10.1} is transformed to
$$
\pp_\rho(\la)
=n^{\frac{|\rho|+\ell(\rho)}2}\,
\prod_{j\ge1}\tf_{j+1}^{\,m_j(\rho)}[\bar\la(\,\cdot\,)]\,
+\,O\left(n^{\frac{|\rho|+\ell(\rho)}2-1}\right),  \tag10.3
$$
Expand $\pp_\rho$ into the sum of its weight homogeneous components:
$$
\pp_\rho=\sum_{j=0}^{|\rho|+\ell(\rho)}F^{(j)}_\rho\,,
\qquad \wt(F^{(j)}_\rho)=|\rho|+\ell(\rho)-j, \tag10.4
$$
so that $F^{(0)}_\rho$ is the top weight component. By Proposition
4.9,
$$
F^{(0)}_\rho =\prod_{j\ge1}(\text{the top weight component of
$\pp_j$})^{m_j(\rho)}\,.
$$
Hence, by Proposition 10.1,
$$
F^{(0)}_\rho=\prod_{j\ge1}\tf_{j+1}^{\,m_j(\rho)}\,. \tag10.5
$$

By virtue of \tht{2.17}, any weight homogeneous element of $\Alg$ is
either symmetric or antisymmetric with respect to ``$\inv$'',
depending on whether its weight is even or odd. It follows from
\tht{4.9} that the element $\pp_\rho$ is either symmetric or
antisymmetric, depending on the parity of the number
$|\rho|+\ell(\rho)$. It follows that, in the expansion \tht{10.4}, we
have $F^{(j)}_\rho=0$ for all odd $j$.

Using this and applying Proposition 2.11, we get from \tht{10.4}
$$
\pp_\rho(\la)=\sum_{j=0}^{[(|\rho|+\ell(\rho))/2]}
n^{\frac{|\rho|+\ell(\rho)}2-j}\,F^{(2j)}_\rho[\bar\la(\,\cdot\,)].
\tag10.6
$$
Set $\om=\bar\la$ and let, as usual, $\si=(\om-|\,\cdot\,|)/2$. By
virtue of the assumption $\la\in\Y(A)$, the support of $\si$ is
contained in $[-A,A]$. It follows that $|\tp_k[\om]|\le 2A^k$ for any
$k\ge2$ (to see this, apply \tht{2.2} and the general estimate
$|\si'(\,\cdot\,)|\le1$). Hence for any element $F\in\Alg$ we get the
estimate
$$
|F[\bar\la(\,\cdot\,)]|\le\operatorname{Const}\,,
\qquad \la\in\Y(A), \tag10.7
$$
where the constant depends only on $A$ and the degree of $F$ as a
polynomial in $\tp_2,\tp_3,\dots$. Applying the estimate \tht{10.7}
to the terms of the expansion \tht{10.6} and taking into account
\tht{10.5} we get the required formula \tht{10.1}. \qed
\enddemo

\enddocument
\end